\documentclass[10pt,oneside,reqno]{amsart}
\usepackage{amsmath,amssymb,amsthm,amsopn,amstext,amsfonts,amsbsy,bbm}
\def\non{\noindent}

\newtheorem{lemma}{Lemma}
\newtheorem{proposition}{Proposition}
\newtheorem{coro}{Corollary}
\newtheorem{remark}{Remark}

\title[Exponential ergodicity  for non-reversible diffusions]{Exponential
  ergodicity and Rayleigh-Schr\"odinger series for  infinite
dimensional diffusions}

\author{ Alejandro F. Ram\' \i rez$^1$
}
\thanks{
$^1$Partially supported by Fondo Nacional de Desarrollo Cient\'\i fico
y Tecnol\'ogico grants 1060738 and 1141094 and
by Iniciativa Cient\'\i fica Milenio NC120062.}

\date{\today}
\email{aramirez@mat.puc.cl}
\address{ Facultad de Matem\'aticas\\
Pontificia Universidad Cat\'olica de Chile\\
Casilla 306-Correo 22, Santiago 6904411, Chile\\
T\'el\'ephone: [56](2)354-5466\\
T\'el\'efax: [56](2)552-5916}

\begin{document}

\begin{abstract} We consider an infinite dimensional diffusion
 on $T^{\mathbb Z^d}$, where $T$ is the circle, defined
by an infinitesimal generator of the form 
$L=\sum_{i\in\mathbb Z^d}\left(\frac{a_i(\eta)}{2}\partial^2_i
+b_i(\eta)\partial_i\right)$, with $\eta\in T^{\mathbb Z^d}$,
 where the coefficients $a_i,b_i$
are  of finite range, bounded with
uniformly bounded second order partial
derivatives and the ellipticity assumption $\inf_{i,\eta}a_i(\eta)>0$
is satisfied. We prove that whenever $\nu$ is an invariant  Gibbs measure
for this diffusion satisfying the logarithmic Sobolev inequality,
then the dynamics is exponentially ergodic in the uniform norm,
and hence $\nu$ is the unique invariant measure. As an
application of this result, we prove that if
 $A=\sum_{i\in\mathbb Z^d}c_i(\eta)\partial_i$,
 and $c_i$ satisfy the condition $\sum_{i\in\mathbb Z^d} \int c_i^2d\nu<\infty$, then there is an $\epsilon_c>0$, such that for
every $\epsilon\in (-\epsilon_c,\epsilon_c)$, the infinite dimensional
diffusion with generator $L_\epsilon=L+\epsilon A$, has a
unique
invariant measure $\nu_\epsilon$ having a Radon-Nikodym
derivative $g_\epsilon$ with respect to $\nu$, which admits
the analytic expansion $g_\epsilon=\sum_{k=0}^\infty \epsilon^k f_k$,
where $f_k\in L_2[\nu]$ are defined through $f_0=1$, 
 $\int f_kd\nu=0$ and the recurrence
equations
$L^*f_{k+1}=A^*f_k$. We give an example where through this expansion
we are able to quantify the effect on the invariant measure of a perturbation triggering
interaction on independent diffusions.
\end{abstract}

\maketitle

\noindent {\footnotesize

}





\newtheorem{teo}{Theorem}
\newtheorem{teof}{Th\'eor\`eme}

\smallskip

\section{ Introduction.}

Let $T$ be the unit circle.  Consider $\Omega:=C([0,\infty );T^{{\mathbb Z}^d})$,
the space of continuous functions from $[0,\infty )$ to
 the set $T^{{\mathbb Z}^d}$,  with the topology of uniform convergence
in compact subsets of $[0,\infty )$. Here
$T^{{\mathbb Z}^d}$ is 
endowed with the product topology. Let $S_t$ be the unique semi-group
on the set $C(T^{{\mathbb Z}^d})$ of continuous real functions on $T^{{\mathbb Z}^d}$
endowed with the uniform norm $\|\cdot\|_\infty$,
associated to the generator which is the closure on $C(T^{{\mathbb Z}^d})$ of the
operator $(L,D_0)$ where
  $L:=\sum_{i\in{{\mathbb Z}^d}}\left(\frac 12a_i (\eta )
\partial_i^2+b_i(\eta)\partial_i
\right)$, with $\partial_i:=\frac{\partial}{\partial\eta_i}$,
 and $D_0$ is the set
of  local functions with continuous second order partial derivatives.
Here, $a
: T^{{\mathbb Z}^d}\to [0,\infty)^{{\mathbb Z}^d}$,
$b: T^{{\mathbb Z}^d}\to \mathbb R^{{\mathbb Z}^d}$ 
 are
Borel-measurable functions which we call sets of {\it coefficients}, $\eta\in T^{{\mathbb Z}^d}$, and
 $a_i, b_i$ and $\eta_i$ are their $i$-th components.
For $i\in\mathbb Z^d$ call its $l_1$ and $l_2$ norms $|i|_1$ and $|i|_2$ respectively.
 We  say that the coefficients $a$ and $b$ are 
 finite range  $R\in{\mathbb Z}^+$ 
if for each
$i\in{\mathbb Z}^d$, $a_i(\eta )$ and $b_i(\eta )$ depend
only on coordinates $\eta_j$ of $\eta$  such that $|j-i|_1\le R$
and
bounded if $sup_{i,\eta}\{a_i,|b_i|\}<\infty$.
 We  say that the
coefficients $a$ and $b$ have uniformly bounded second order
partial derivatives if
 $\sup_{i,j,k,\eta}\left\{\left|\frac{\partial^2a_i}{\partial\eta_j\partial\eta_k}\right|,
\left|\frac{\partial^2b_i}{\partial\eta_j\partial\eta_k}\right|\right\}<\infty$.
It is easy to check that whenever $a$ and $b$ have uniformly bounded second
order partial derivatives, they also have uniformly bounded first order derivatives,
so that
 $\sup_{i,j,\eta}\left\{\left|\frac{\partial a_i}{\partial\eta_j}\right|,
\left|\frac{\partial b_i}{\partial\eta_j}\right|\right\}<\infty$.
Throughout the rest of
this paper, we will always consider coefficients which are 
 finite range, bounded and with uniformly bounded
second order partial derivatives, and we will call them
just {\it coefficients}. 
It is a standard fact (which for completeness will be shown
in Section \ref{section2}), that the operator $(L,D_0)$ defined above,
with coefficients  $a$ and
$b$  is closable and of Hille-Yosida type. We will also see
that whenever we also have that

\begin{equation}
\label{uel}
a:=\inf_{i,\eta}a_i(\eta)>0,
\end{equation}
its closure, which we will
denote 

\begin{equation}
\label{closure}
(L,D(a,b)),
\end{equation} is an infinitesimal generator and defines a
Markov semi-group $\{S_t:t\ge 0\}$  on the space 
$C(T^{{\mathbb Z}^d})$ corresponding to a diffusion
process on $\Omega$
(see  \cite{stroock} where an alternative construction is performed
solving the martingale problem for $(L,D_0)$). 
Such a  process will be called a
 {\it infinite dimensional diffusion with coefficients $a$ and $b$}
or just {\it infinite dimensional diffusion}.
A coefficient $a$ satisfying (\ref{uel}) will be called
a {\it uniformly elliptic} coefficient.

Given a probability measure $\nu$ defined on $T^{{\mathbb Z}^d}$ endowed
with its Borel $\sigma$-algebra,
throughout the paper we will use the notation $\langle f\rangle_\nu:=
\int f d\nu$ for $f$ integrable defined in $C(T^{\mathbb Z^d})$. 
For each $p\ge 1$, 
 we denote by $L_p[\nu]$ the
Banach space of functions $f: T^{\mathbb Z^d}\to\mathbb C$ 
with norm $\|f\|_{p,\nu}:=\left(\int |f|d\nu)\right)^{1/p}$.
We adopt the convention that $L_2[\nu]$ is the space
of complex valued square integrable functions
and given two functions $f,g\in L_2[\nu]$, we denote their
inner product by $(f,g)_\nu:=\int \bar fgd\nu$. 
Furthermore, for every $t\ge 0$, we define $\nu S_t$ as the unique measure such that
$\int S_t fd\nu=\int fd (\nu S_t)$ for every continuous function $f$
on $T^{{\mathbb Z}^d}$.
We will call this measure
an {\it invariant measure} for the infinite dimensional diffusion, if for every $t\ge 0$
one has that $\nu S_t=\nu$.

Note that due to compactness of the state space,
an infinite dimensional diffusion has always at least
one invariant measue. Nevertheless, few general results exist providing sufficient conditions for
the existence of a unique invariant measure, or for the exponential
ergodicity of infinite dimensional diffusions, specially out
of the subclass of reversible processes. Furthermore, given that
in general it is difficult to explicitly describe explicitely the invariant
measures, it is natural to wonder what is the breadth
of for example, the classical theory of
analytic perturbations for the invariant 
measures. In this paper we address these issues, within
the context of Gibbs probability measures satisfying the logarithmic Sobolev
inequality and not satisfying any kind of reversibility assumptions.
Our results are  fundamental  settling down
uniqueness and analyticity issues under general conditions.

We recall the definition of the logarithmic Sobolev inequality which
will be assumed throughout this article.

\smallskip
\noindent{\bf Condition (LSI)}. {\it
We say that a probability measure $\nu$ on $T^{{\mathbb Z}^d}$ endowed
with its Borel $\sigma$-field, satisfies the logarithmic Sobolev inequality 
with respect to the Laplacian operator if there is a constant $\gamma>0$ such
that for every non negative function $f\in D_0$  it is true that

\begin{equation}
\label{lsi}
\left\langle f^2\ln\frac{f}{\sqrt{\langle f^2\rangle_\nu]}}\right\rangle_\nu
\le \gamma
\left\langle\sum_{i\in{{\mathbb Z}^d}}  \left(\partial_i f\right)^2\right\rangle_\nu.
\end{equation}}

\smallskip

Let us also recall the definition of a Gibbs measure. A potential on $T^{\mathbb Z^d}$
is a collection $\mathcal J$ of functions $\{J_F:F\subset\mathbb Z^d, F\ {\rm finite}\}$,
such that for each $F$ the function $J_F$ has continuous second order
partial derivatives and is invariant under permutations
of the indices of $F$. We will also assume that the collection
$\mathcal J$ is {\it finite range}: there exists an $L$ such that
$i,j\in F$ with $|i-j|_1>L$ implies that $J_F=0$. Define now the
energy of the subset $G\subset\mathbb Z^d$ as

$$
H_G(\eta):=\sum_{F:  G\subset F} J_F(\eta).
$$
We now say that a probability measure $\mu$ on $T^{\mathbb Z^d}$ is
a {\it Gibbs measure with potential $\mathcal J$} if for all
$G\subset\mathbb Z^d$ the regular probability distribution
of $\mu$ given $\sigma(\eta_l:l\notin G)$ admits a
density $u_G(\{\eta_k:k\in G\})$ given by

$$
u_G(\{\eta_k:k\in G\})=\frac{e^{-H_G(\eta)}}{\int e^{-H_G(\zeta)}d\prod_{k\in G}\zeta_k}.
$$

\smallskip

On the other hand, we say that $\nu$ satisfies
the spectral gap inequality with spectral gap constant $g$ with respect to the
Laplacian operator (\ref{lsi}) if for every $f\in D_0$
one has that

\begin{equation}
\label{thespectral}
g||f-\langle f\rangle_\nu||^2_{2,\nu}\le \left\langle\sum_{i\in{\mathbb Z}^d}\
(\partial_i f )^2\right\rangle_\nu
\end{equation}
(see \cite{gz}). It is a standard fact, which we will subsequently use,
that whenever $\nu$ satisfies the logarithmic Sobolev inequality with
constant $\gamma$ it
necessarily satisfies the spectral gap  inequality with spectral gap constant $g=1/\gamma$.

 The first result of this paper provides a sufficient
condition for
exponential ergodicity of infinite dimensional diffusions.
Let us define  for each $\theta>0$, the
triple semi-norm for  $f\in D_0$ as

\begin{equation}
\label{triple-norm}
|||f|||_\theta:=\left(\sum_{i,j\in\mathbb Z^d}e^{\theta (|i|_2+|j|_2)}\sup_{\eta}\left(
\partial_i\partial_j f\right)^2\right)^{1/2}
\end{equation}
(this type of semi-norm was already used in \cite{ramirez} to control
the dependence of functions evolved according to the semi-group of
infinite dimensional diffusions). Remark that since
$f\in D_0$, the sum in (\ref{triple-norm}) is finite.

\smallskip
\begin{teo} 
\label{theorem2}  Consider an infinite dimensional diffusion with
semi-group $\{S_t:t\ge 0\}$ and with an invariant Gibbs
measure $\nu$ 
which
satisfies {\bf (LSI)}. Let $\theta>0$. Then, there exist positive constants $k_\theta$ and
$K_\theta$, depending only on $\theta$ and on the coefficients
of the diffusion, such that for
any  function $f\in D_0$ we have,

$$\sup_{\eta\in T^{{\mathbb Z}^d}}\left|S_t f(\eta )-\langle f\rangle_\nu \right|\le 
K_\theta \max\{|||f|||_\theta,1\} e^{-k_\theta t}.$$
\end{teo}

\smallskip
\noindent Theorem \ref{theorem2} is an improvement of a result of
Zegarlinski \cite{z} for reversible processes.
 As a corollary of Theorem \ref{theorem2},
 we obtain the following considerable improvement
of Theorem 1 of  \cite{r1}. 

\smallskip

\begin{coro}
\label{cor}
Consider an  infinite dimensional 
diffusion with an invariant Gibbs measure $\nu$ which satisfies
{\bf (LSI)}. Then, $\nu$  is unique and 
for every probability measure $\mu$ on
$T^{{\mathbb Z}^d}$, one has that $\lim_{t\to\infty}\mu S_t=\nu$ weakly.
\end{coro}

\smallskip


\noindent Theorem \ref{theorem1}, which we formulate below,  establishes
both a stability result for the uniqueness of invariant Gibbs measures
of infinite dimensional diffusions
satisfying the logarithmic Sobolev inequality under
a certain class of perturbations, and
 the existence of a Rayleigh-Schr\"odinger
series with a positive radius of convergence around 
the invariant measure of the unperturbed diffusion.
To state Theorem \ref{theorem1}, we need to introduce the
following regularity condition on coefficients.

\smallskip
\noindent{\bf Condition (R)}. {\it
Let $\mu$ be a probability measure defined on $T^{\mathbb Z^d}$. We say that
a set of  coefficients $c=\{c_i:i\in
\mathbb Z^d\}$ satisfies the regularity  condition {\bf (R)} with
respect to $\mu$, if

\begin{equation}
\label{ce0}
C_0:=\sqrt{\sum_i \int c_i^2 d\mu}<\infty.
\end{equation}}
\medskip

\non Note that a sufficient condition for condition {\bf (R)} to
be satisfied with respect to any probability measure $\mu$  is that

$$
\sup_{\eta\in T^{\mathbb Z^d}}\sum_i c_i^2<\infty.
$$
Given a set of coefficients $c$ satisfying condition {\bf (R)} with
respect to $\mu$,
any operator of the form 

\begin{equation}
\label{dfo}
A=\sum_{i\in\mathbb Z^d}c_i(\eta)\partial_i,
\end{equation}
will be called a {\it diagonal first order operator} satisfying condition {\bf (R)} with respect to $\mu$.

It will be
shown,  that if $L_0$ is the generator of an
  infinite dimensional diffusion with invariant measure $\nu$
satisfying {\bf (LSI)}
and $A$ a diagonal first order operator satisfying condition {\bf (R)} with respect to $\nu$,
 then the closure of the
operators $(L_0,D_0)$ and $(L_0+A,D_0)$, on $L_2[\nu]$ have the same domain,
which for the moment we will call $\bar D_0$.
Furthermore, given an operator $(T,D(T))$ on $L_2[\nu]$, we
denote by $(T^*,D(T^*))$ its adjoint. Hence,
the adjoints of $(L_0,\bar D_0)$ and $(A,\bar D_0)$ have
the same domain, which we will denote by $\bar D_0^*$.
We will now see that for every function $f\in\bar D_0^*$, the equation

\begin{equation}
\label{weak}
L^*_0g=-A^*f,
\end{equation}
has a solution $g\in \bar D_0^*$. 
 Here, equation (\ref{weak}) should
 be interpreted in the weak sense as $(g,L_0\phi)_\nu=-(f,A\phi)_\nu$,
for every $\phi\in D_0$.
 To see that (\ref{weak}) has a solution, it
is enough to prove that $A^*f$ is orthogonal to the kernel
of $L^*_0$. But it will be shown that $0$ is a simple
eigenvalue of $L^*_0$ and that the eigenfunctions of $0$ are the
constant ones. So it is enough to prove that $A^*f$ is orthogonal
to constants. Now,

$$
(1,A^*f)_\nu=(A1,f)_\nu=0,
$$
which proves the statement.
The fact that the eigenvalue $0$ of $L^*_0$ is  simple
 implies that the solution $g$ of the equation
(\ref{weak}) is unique and that $\langle g\rangle_\nu=0$.
We will denote by $M_2[\nu]$  the set of probability measures
on $T^{{\mathbb Z}^d}$ which have a Radon-Nikodym derivative with
respect to $\nu$ which is square integrable in $L_2[\nu]$.
Let

\begin{equation}
\label{epsc}
\epsilon_c:=\frac{a}{C_0\sqrt{\gamma}}.
\end{equation}

\smallskip
\begin{teo}
\label{theorem1} Consider
an  infinite dimensional diffusion
with infinitesimal generator $L_0$ and with
an invariant Gibbs  measure  $\nu$ 
satisfying {\bf (LSI)}.
  Let $A$ be a diagonal first order operator
satisfying condition {\bf (R)} with respect to $\nu$.
Then,  for each $\epsilon\in (-\epsilon_c,
\epsilon_c)$, 
 the infinite dimensional diffusion
with generator $L_\epsilon:=L_0+\epsilon A$ has a unique invariant
 measure $\nu_\epsilon$ in $M_2[\nu]$ which has a  Radon-Nikodym derivative
 $g_\epsilon$  with respect to $\nu$ with the  following  expansion
in $L_2[\nu]$,

\begin{equation}
\label{gexpansion}
g_\epsilon=\sum_{k=0}^\infty \epsilon^k f_{k},
\end{equation}
where $\{f_k:k\ge 0\}$ is the unique sequence of functions in $L_2[\nu]$ 
defined by $f_0:=1$,
the conditions $\langle f_k\rangle_\nu=0$ for $k\ge 1$,
 and the recurrence relations

$$
L_0^*f_{k+1}=-A^*f_k.
$$
Furthermore, there exists a constant $C$ such that for every $k\ge 1$ one
has that $\|f_k\|_{2,\nu}\le C\epsilon_c^{-k}$.

\end{teo}
\smallskip

\non Theorem \ref{theorem1} does not require the unperturbed generator to be
reversible with respect to the invariant measure. In \cite{ko}, within the
context of systems which satisfy the Einstein relation, a similar
expansion was derived for interacting particle systems, under the assumption
that the unperturbed generator is reversible. It should also be
noted that in many interesting situations, for example within the context
of random walks in random environments, one cannot expect in general an analytic
expansion of the invariant measure (see \cite{cr15}).

 On the other hand, note
that the diagonal first order operator $A$ is not bounded (hence not compact) in $L_2[\nu ]$.
Furthermore, we do not know if condition ${\bf (R)}$ in
Theorem \ref{theorem1} could be relaxed, but it is possible that
any essential generalization of such a condition would
break up the analiticity of the expansion (\ref{gexpansion}).

A crucial ingredient in the proofs of Theorems \ref{theorem2} and \ref{theorem1},  is the
observation that the Dirichlet form of the generator of a diffusion
having an invariant measure that satisfies the logarithmic Sobolev inequality,
coincides with the Dirichlet form of the symmetrization of the generator.
For Theorem \ref{theorem2}, this together with the spectral gap, implies an exponentially fast
convergence to the equilibrium measure in the $L_2[\nu]$ norm. 
One can then get exponentially fast convergence in the supremum norm
through the following three additional ingredients: comparisons
between a truncated version of the dynamics and the full dynamics, Gross lemma
and uniform norm estimates on the marginal distribution of the process which
are obtained using Girsanov theorem. Some of these ingredients, as for example
Gross lemma, are not new, while others as the uniform norm estimates for marginal distributions, to the knowledge of the author is novel, requiring a careful
use of Girsanov's theorem. Nevertheless, somehow surprisingly, they
have not been combined in this way before to provide a non-reversible result.
The proof of Theorem \ref{theorem1},
uses the uniqueness result of Theorem \ref{theorem2}
for invariant measures satisfying the logarithmic Sobolev
inequality (Corollary \ref{cor}) and uses the machinery of
analytic perturbation theory for operators which have a relatively bounded
perturbation on Banach spaces. In particular, the uniqueness of
the perturbed invariant measure in $M_2[\nu ]$ of Theorem \ref{theorem1}
will be derived using Corollary \ref{cor}.

In Theorem \ref{construction} of the next section, we construct
 infinite dimensional diffusions via Hille-Yosida theorem.
 In section 3, we derive some important
consequences of the property that a diffusion has an invariant
measure satisfying the  spectral
gap inequality.  In section 4, we use the results of section 3, to prove
Theorem \ref{theorem2}. Theorem \ref{theorem1} is proved in section 5.
In section 6,  Theorem \ref{theorem1} is illustrated
applying  giving a structural expansion
 for the invariante measure of weakly perturbed
independent diffusions.

\smallskip

\section{Notation and preliminaries}
\label{section2}
 
In this section,  we will provide a construction of
 infinite dimensional diffusions 
via the Hille-Yosida theorem. In fact, we will show that every
operator $(L,D_0)$ of an infinite dimensional diffusion
is closable, and that its closure defines a Markov generator
of a Markov semi-group $\{S_t:t\ge 0\}$ on $C(T^{\mathbb Z^d})$. Classical
results in the theory of Markov processes then imply
that there exists a Markov process $\{P_\eta:\eta\in T^{\mathbb Z^d}\}$
such that $S_tf(\eta)=E_\eta[f(\eta_t)]$ for all $f\in C(T^{\mathbb Z^d})$,
$\eta\in T^{\mathbb Z^d}$ and $t\ge 0$ (see for example Theorem 1.5
of Liggett \cite{liggett} or Dynkin \cite{dynkin}).
Infinite dimensional diffusions have already been constructed via the martingale
problem by Holley and Stroock \cite{stroock}, and we do not claim
that our presentation  is particularly original (see also
\cite{zak} for a more recent reference). Nevertheless,
for the sake of completeness, and lacking an appropriate reference
for the construction  via the Hille-Yosida theorem,
we have included the details here. Throughout, we
will denote by $\mathbb I$ the identity operator on $C(T^{\mathbb Z^d})$.
For $r>0$, we will define the box $\Lambda_r:=[-r,r]^d\cap{\mathbb Z}^d$.
Given a closable operator $(A,D)$ on $C(T^{\mathbb Z^d})$,
we will denote by $(A,\tilde D)$ its closure. 

\smallskip

We will also need to introduce the notion of truncated operator.
For each subset $\Lambda\subset\mathbb Z^d$, define 
$\mathcal F_{\Lambda}$ as the $\sigma$-algebra of
Borel subsets of $T^{\mathbb Z^d}$ generated by the coordinates
in $\Lambda$.
Let us fix a natural $n\ge 1$ and a probability
measure $\mu$ on $T^{\mathbb Z^d}$.
Given two sets of coefficients
$a$ and $b$, and given the corresponding operator of the form $(L,D_0)$ with
$L=\sum_{i\in\mathbb Z^d}\left(\frac{1}{2}a_i(\eta)\partial_i^2
+b_i(\eta)\partial_i\right)$,  
we  call the pair $(L_n,D_n)$
its {\it truncation integrated with respect
to the measure $\mu$}, where
 $D_n$ is the set of
functions in $D_0$ 
which are $\mathcal F_{\Lambda_n}$ measurable (or equivalently, the functions depending
only on the coordinates in $\Lambda_n$ with continuous second order partial
derivatives) 
and where

\begin{equation}
\label{tg}
L_n: =\sum_{i\in\Lambda_n}\left(
\frac{1}{2}\bar a_i\partial^2_i +\bar b_i\partial_i \right),
\end{equation}
with 

\begin{equation}
\label{tcoef}
\bar a_i(\eta):=E_\mu[ a_i(\eta )|\mathcal F_{\Lambda_n}]
\quad {\rm and}\quad 
\bar b_i(\eta):=E_\mu[ b_i(\eta )|\mathcal F_{\Lambda_n}],
\end{equation} the
conditional expectation of $a_i$ and $b_i$ with respect to $\mu$
 given the $\sigma$-algebra ${\mathcal F}_{\Lambda_n}$. We
call $\bar a$ and $\bar b$ the corresponding coefficients. 
Given $\eta\in T^{{\mathbb Z}^d}$,  we define $\eta^n$, called the {\it
  configuration $\eta$ truncated at
scale $n$}  by 

\begin{equation}
\label{eta-trunc}
\eta^n(x):=
\begin{cases}
\eta(x)\quad &{\rm if}\quad |x|\le n\\
1\quad &{\rm otherwise}.
\end{cases}
\end{equation}
For a function $f\in D_0$ and $i\in\mathbb Z^d$, we define

$$
\Delta_f(i):=\sup_\eta|\partial_i^2f|.
$$
This is a measure of the dependance of $f$ on its $i$-th coordinate (see
\cite{ramirez}, where this quantity was already introduced).

It is a standard fact from the theory
of partial differential equations that for each $n$, the truncated
operator $(L_n,D_n)$ is closable in $C(T^{\Lambda_n})$, and that its closure,
which we will call

\begin{equation}
\label{tgen}
(L_n,D_n(a,b)),
\end{equation}
  is a Markov generator (see \cite{liggett}
for the definition of Markov generator) and hence defines
a semi-group (see \cite{pinsky} or Chapter 12 of \cite{varadhan})

\begin{equation}
\label{semin}
\{S_t^n:t\ge 0\}\quad {\rm on\ the\ space}\quad C(T^{\Lambda_n}).
\end{equation}
We will require the following version of Theorem 2 of \cite{ramirez}.

%

\smallskip

\begin{teo}
\label{theorem1.3} Consider coefficients $a$ and $b$, with
$a$ uniformly eliptic, and for each $n\ge 1$
its truncated operators
$(L_n,D_n)$  with respect to Lebesgue measure. Let
$1\le M\le n$. Then, for each $\vartheta>0$, $g\in D_M$  and $i\in\mathbb Z^d$, there exist constants $A$ and $\beta$
(depending only on $M$, $\vartheta$  and the coefficients
of the infinite dimensional diffusion) such that

$$
\Delta_{ S_t^n g}(i)\le A |||g|||_\vartheta e^{\beta t}e^{-\vartheta |i|_2},
$$
for all $t\ge 0$.
\end{teo}

\smallskip

\noindent The proof of \ref{theorem1.3} is completely analogous
to that of Theorem 2 in \cite{ramirez}, but for completeness
we present it in the Appendix \ref{appendix}.
With the control provided by Theorem \ref{theorem1.3} we can
now prove the $(L,D_0)$ is closable and that its closure is
a Markov generator (see \cite{liggett} for the definition of
Markov generator).

\smallskip
\begin{teo}
\label{construction} 
Any operator $(L,D_0)$ with coefficients $a$ and $b$, with
$a$ uniformly elliptic, is closable in $C(T^{{\mathbb Z}^d})$,
and its closure is the Markov generator of a Markov semi-group
on $C(T^{{\mathbb Z}^d})$.
\end{teo}
\begin{proof} Let us first note that by Propositions 2.2 and 2.5 of
\cite{liggett}, to show that $(L,D_0)$ is closable it is enough
to prove that for every $f\in D_0$
and $\zeta\in T^{\mathbb Z^d}$ such that
 $\inf_{\eta\in T^{\mathbb Z^d}}f(\eta)= f(\zeta)$, one has that
$Lf(\zeta)\ge 0$. Now, note that for each $i\in\mathbb Z^d$
it is true that $\frac{\partial f}{\partial\eta_i}(\zeta)=0$
and that $\frac{\partial^2 f}{\partial^2\eta_i}(\zeta)\ge 0$,
so the statement is proved.
Next, to prove that the closure of $(L,D_0)$
is a Markov generator, it is enough to show that for every $\lambda>0$ small enough,
the range of the operator $\mathbb I-\lambda L$ is dense in
$C(T^{\mathbb Z^d})$.  For this purpose, consider for each natural $n$
the truncated version $L_n$ of $L$ integrated with respect to
Lebesgue measure. Let $g\in C(T^{\mathbb Z^d})$ and $\epsilon>0$.
Now choose $M$ such that if  $g_M(\eta):=g(\eta^M)$, we have
that

\begin{equation}
\label{gem}
||g-g_M||_\infty\le\frac{\epsilon}{3}.
\end{equation}
Then choose $h_M\in D_M$ so that

\begin{equation}
\label{gemh}
||g_M-h_M||_\infty\le\frac{\epsilon}{3}.
\end{equation}
Note that by the fact that $h_M$ and the coefficients
of the operator $L_n$ are differentiable with
derivatives which are H\"older continuous of positive parameter,
by classical regularity theory for elliptic partial
differential equations
(see for example Theorems 6.14 and 6.17 of \cite{gt} which can be
easily adapted to the periodic domain $T^{\Lambda_n}$),
for each  $n\ge M$, there exists a function
$f_n\in D_n$ such that 

\begin{equation}
\label{resolv}
(\mathbb I-\lambda L_n)f_n=h_M.
\end{equation}
Furthermore, we have the following representation for the
solution of (\ref{resolv}),

\begin{equation}
\label{identt}
f_n=-\frac{1}{\lambda}\int_0^\infty e^{-\frac{t}{\lambda}}S_t^n h_Mdt,
\end{equation}
where $\{S_t^n:t\ge 0\}$ is the semi-group (\ref{semin}) of the corresponding
finite dimensional diffusion (see for example Theorem 1 of chapter 12 of \cite{varadhan}).
On the other hand, it is easy to check that

\begin{equation}
\label{lfn}
||Lf_n-L_nf_n||_\infty\le C\sum_{i\notin\Lambda_{n-R}}\Delta_{f_n}(i),
\end{equation}
where for each $i\in\mathbb Z^d$, $\Delta_{f_n}(i)<\infty$ since $f\in D_n$.
By Theorem \ref{theorem1.3}, the identity (\ref{identt}) and
the fact that $h_M\in D_n$, it follows that for given $\vartheta>0$,
we can find constants $\beta>0$ and $A$ depending only
on $\vartheta$ and $M$, such that for every $\lambda<\frac{1}{\beta}$
one has that

$$
\Delta_{f_n}(i)\le\frac{1}{\lambda}\int_0^\infty e^{-\frac{t}{\lambda}}\Delta_{S_t^nh_M}(i)dt
\le \frac{A}{\lambda}|||h_M|||_\vartheta e^{-\vartheta |i|_2}\int_0^\infty e^{-\left(\frac{1}{\lambda}-\beta\right)t}dt 
=
|||h_M|||_\vartheta \frac{A}{1-\lambda\beta} e^{-\vartheta |i|_2}.
$$
Substituting this back into inequality (\ref{lfn}), we get that

$$
||Lf_n-L_nf_n||_\infty\le C|||h_M|||_\vartheta \frac{A}{1-\lambda\beta} \sum_{i\notin\Lambda_{n-R}}e^{-\vartheta |i|_2},
$$
which implies that for $n$ large enough one has that $||Lf_n-L_nf_n||_\infty\le\frac{\epsilon}{3}$.
Combining this with inequalities (\ref{gem}) and (\ref{gemh}), we conclude that there exists a $\lambda_0>0$ such that
for each $\lambda<\lambda_0$ we have that for $n$ large enough

$$
||(\mathbb I-\lambda L)f_n-g||_\infty\le\epsilon.
$$
This proves that for $\lambda>0$ small enough, the range of $\mathbb I-\lambda L$
is dense in $C(T^{\mathbb Z^d})$.
\end{proof}

\noindent  Recall that given an operator $(L,D_0)$ with coefficients
$a$ and $b$, its closure is denoted by $(L,D(a,b))$ [c.f. (\ref{closure})].
By Theorem \ref{construction}, when $a$ is uniformly elliptic, $(L,D(a,b))$ is a Markov generator
of an infinite dimensional difusion.
If $\mu$ is an invariant measure for this diffusion and $\{S_t:t\ge 0\}$ is
its Markov semi-group, then  by Jensen inequality for every function
$f\in C(T^{{\mathbb Z}^d})$ one has that,

\begin{equation}
\nonumber
||S_t f||_{2,\mu}\le ||f||_{2,\mu}.
\end{equation}
It is hence possible to continuously extend the Markov semi-group $S_t$ as a Markov semi-group
to $L_2[\mu]$ and with an infinitesimal generator
which is  an extension is the closure in $L_2[\mu]$
of the  generator $(L, D_0)$
(see for example Proposition 4.1 of \cite{liggett}). We will denote it $(L,\bar D(a,b))$. 
For the truncated generator $L_n$ of $L$, we will denote
the closure of $(L_n,D_n)$ in $L_2[\mu]$ as $(L_n,\bar D_n(a,b))$
and its Markov semi-group on $L^2[\mu]$ by $\{S_t^n:t\ge 0\}$
(we thus make a slight abuse of notation not distingushing between this
semi-group and the one in $C(T^{\Lambda_n})$ defined in (\ref{semin})).
In the sequel, we will make no notational distinction between
the semi-group on $C(T^{{\mathbb Z}^d})$ or on $L_2[\mu]$.
Given $\Lambda\subset\mathbb Z^d$, we will denote by $\mu_\Lambda$ the restriction
of $\mu$ to the $\sigma$-algebra ${\mathcal F}_\Lambda$.

\smallskip

\section{Symmetrization of the infinitesimal generator.} 
Here, we will show that an infinite dimensional
diffusion has a Dirichlet form which is equal
to the Dirichlet form of its symmetrized generator
with respect to a given invariant measure.
If the invariant measure satisfies the spectral
gap inequality [c.f. (\ref{thespectral})], this in turn   implies an exponential convergence result
 to equilibrium in the corresponding
$L_2$ norm. We will first need the following lemma.

\smallskip

\begin{lemma}
\label{lemma-new} Consider an infinite dimensional diffusion
with infinitesimal generator  $(L, D(a,b))$ and a probability measure
$\mu$ on $T^{\mathbb Z^d}$. Then, for
every real function $f\in D_0$ one has that

\begin{equation}
\label{key-identity}
(f,Lf)_\mu=\frac{1}{2}\int L f^2d\mu-\frac{1}{2}\sum_{i\in\mathbb Z^d}\int a_i(\partial_i f)^2
d\mu.
\end{equation}
\end{lemma}

\smallskip

\begin{remark} The second term of the right-hand side of (\ref{key-identity})
is the negative of the integral of the so called  {\it ``carr\'e du champ''} of the diffusion with
generator $L$.
\end{remark}

\smallskip
\begin{proof}[Proof of Lemma \ref{lemma-new}] Let us choose $n$ natural so that the support
of $f$ is in $\Lambda_n$. It follows that $fLf$ is a local
function with support in $\Lambda_{n+R}$. Therefore, to prove
(\ref{key-identity}) we can assume that $\mu$ is in fact a probability
measure on $C(T^{\Lambda_{n+R}})$. On the other hand,
if $\mu$ is a probability measure with a smooth density $u$
with respect to the Lebesgue measure $m$ on $T^{\mathbb Z^d}$, we have that

$$
(f,Lf)_\mu=\int fLf udm,
$$
and then
 (\ref{key-identity})
can be easily deduced through integration by parts.
We can now
conclude the proof observing that given an arbitrary probability
measure $\mu$ defined on $C(T^{\mathbb Z^d})$, we can always construct
a sequence of measures which have smooth densities with respect to
Lebesgue measure, converging weakly to $\mu$.
\end{proof}

\smallskip

\noindent We have now as  consequence of Lemma \ref{lemma-new}
the following proposition.
\smallskip

\begin{proposition}
\label{prop1}
Consider
an infinite dimensional diffusion  with  infinitesimal generator $(L,D(a,b))$.
Let $\nu$ be an invariant measure of this infinite dimensional
diffusion and $(L,\bar D(a,b))$ be its infinitesimal
generator on $L_2[\nu]$.  Then, the following are satisfied.

\begin{itemize}

\item[(i)] For every real $f\in D_0$ we have that

$$
(f,Lf)_\nu=-\frac{1}{2}\sum_{i\in\mathbb Z^d}\int a_i (\partial_i f)^2d\nu.
$$

\item[(ii)] Assume that the invariant measure $\nu$ satisfies the
spectral gap inequality with constant $g$. Then, for every real $f\in \bar D(a,b)$ we have
that 
$$
\|f-\langle f\rangle_\nu\|_{2,\nu}^2\le-\frac{2}{ag}(f,Lf)_\nu.
$$

\end{itemize}
\end{proposition}
\smallskip

\begin{proof}
{\it Part $(i)$.}
This is a direct consequence of Lemma \ref{lemma-new} and
the assumption that $\nu$ is an invariant measure.

\medskip

\noindent {\it Part $(ii)$.} Combining the spectral
gap inequality with part $(i)$, and using
the fact that the diffusion is uniformly elliptic, we
conclude that for every $f\in D_0$,

$$
\|f-\langle f\rangle_\nu\|_{2,\nu}^2\le-\frac{2}{ag}(f,Lf)_\nu.
$$
Using the fact that $D_0$ is dense in $\bar D(a,b)$ we obtain
part $(ii)$.

\end{proof}

\non We continue with the following important corollary of part $(i)$ of
the above proposition.
Recall that 
$a$ is the lower bound defined in (\ref{uel}).

\smallskip

\begin{coro}
\label{coro0}
Consider
an infinite dimensional diffusion  with coefficients $a$ and $b$
and infinitesimal generator $(L,D(a,b))$.
Let $\nu$ be an invariant measure of this infinite dimensional
diffusion and $(L,\bar D(a,b))$ be its infinitesimal
generator on $L_2[\nu]$. 
 Let $A=\sum c_i\partial_i$ be a diagonal
first order operator satisfying condition {\bf (R)} with respect to $\nu$.
Then, for every $\lambda$ positive we have that 
for every $f\in D_0$ the following inequality is satisfied,

\begin{equation}
\label{rs-1}
||Af||_{2,\nu}\le 
C_0\sqrt{\frac{2}{a}}\frac{1}{\lambda} \|Lf\|_{2,\nu}+
C_0\sqrt{\frac{2}{a}}\lambda
||f||_{2,\nu},
\end{equation}
where
$C_0$ is defined in (\ref{ce0}).
In particular, the operator $(A,D_0)$ is relatively bounded with respect to
  $(L,D_0)$ in $L_2[\nu]$, with $L$-bound $0$. Hence the
closure of $(L+A,D_0)$ has the same domain $\bar D(a,b)$ as $L$.
\end{coro}
\begin{proof} 
 Let $f\in D_0$. Let us first assume that $f$ is real valued. Note that

\begin{eqnarray}
\nonumber
&||A f||_{2,\nu}^2=\displaystyle\int
\left(
\sum_{i} c_i\partial_i  f
\right)^2
d\nu\le
C_0^2\displaystyle\sum_i\int(\partial_i f)^2d\nu\\
\nonumber
&\le\displaystyle
-2\frac{C_0^2}{a}  \left(f, Lf\right)_\nu
\le 2\frac{C_0^2}{a}  
\left\| f\right\|_{2,\nu} \left\| Lf\right\|_{2,\nu}.
\end{eqnarray}
Here, we have used part $(i)$ of Proposition \ref{prop1} in the second
 inequality. By Cauchy-Schwartz's inequality it
follows that  for every positive $\lambda$ one has that
$||A f||_{2,\nu}^2\le
2\frac{C_0^2}{a}\left(\lambda^{-2}\left\|
    Lf\right\|^2_{2,\nu}
+\lambda^2\left\| f\right\|^2_{2,\nu} 
\right),$
which proves inequality (\ref{rs-1}). 
By taking $\lambda$ arbitrarily large, we see that $A$ has $L$-bound $0$.
Finally, by Theorem 1.1
of page 190 of Kato \cite{kato}, we can see
that the closure of the operators $(L,D_0)$ and $(L+A,D_0)$ in $L_2[\nu]$
have the same domain. To prove the result for an arbitrary complex
valued $f\in D_0$, it is
enough to use $||f||_{2,\nu}=||Re f||_{2,\nu}+||Im f||_{2,\nu}$
and the inequality just proved for real valued functions.
\end{proof}

\non We also have the following corollary of Proposition \ref{prop1}.

\begin{coro}
\label{lemma1} Consider an infinite dimensional diffusion with coefficients
$a$ and $b$ and infinitesimal generator $(L,D(a,b))$. Let $\nu$
be an invariant measure of this infinite dimensional diffusion
and  $(L,\bar D(a,b))$ be its infinitesimal generator on $L^2[\nu]$.
Assume that $\nu$ satisfies  the spectral gap inequality
with constant $g$ with respect to the Laplacian operator.
 Then, for every $f\in L_2[\nu]$ it is true that,

\begin{equation}
\label{one}
||S_t f-\langle f\rangle_\nu||_{2,\nu}\le e^{-ga t}
||f-\langle f\rangle_\nu||_{2,\nu}.
\end{equation}
\end{coro}

\begin{proof}
Note that for every $f\in L^2[\nu]$  and $t>0$

\begin{eqnarray*}
\frac{d}{dt}||S_t f-\langle f \rangle_\nu ||_{2,\nu}^2
&=&
2\int (S_t f-\langle f\rangle_\nu)L (S_t f-\langle f\rangle_\nu)d\nu\\
&\le &- ag||S_t f-\langle f\rangle_\nu||_{2,\nu}^2,
\end{eqnarray*}
where in the inequality we have used part $(ii)$ of Proposition \ref{prop1},
and the fact that for every $t>0$, $S_tf\in \bar D(a,b)$. From this
inequality we immediately deduce (\ref{one}).
\end{proof}
\smallskip

\section{Proof of Theorem \ref{theorem2}.} 
 The basis of the proof of Theorem \ref{theorem2} is
 Corollary \ref{lemma1} of the
previous section, a truncation estimate, Gross lemma and
a uniform estimate on marginal distributions. Throughout, $\gamma$ will denote the constant appearing in the 
logarithmic Sobolev inequality  (\ref{lsi}). Furthermore, we will
adopt the convention that given any sequence $\{y_n\}$, and positive real number
 $x$,
$y_x:=y_{\lfloor x\rfloor}$.

\medskip

\begin{teo}
\label{lemma2}
 {\bf [Truncation estimate]}. Let $\mu$ be any probability measure
defined on $T^{\mathbb Z^d}$. Consider an infinite dimensional diffusion
with semi-group $\{S_t:t\ge 0\}$ and its truncated semi-group at scale $n$
with respect to $\mu$,
 $\{S_t^n:t\ge 0\}$.
 Let $\theta>0$. Then, for every $\delta>0$ there exist 
constants $c_\theta>0$ and $C_\theta>0$, depending only
on $\theta, \delta$ and the coefficients of the diffusion,
such that for all $f\in D_0$ the following statements are satisfied.

\begin{itemize}

\item[(i)] For every $n\ge c_\theta t$,

\begin{equation}
\nonumber
\sup_{\eta\in T^{{\mathbb Z}^d},u:0\le u\le t} |S_u f(\eta)-S^n_u f(\eta)|\le
 C_\theta |||f|||_\theta e^{-\delta t}.
\end{equation}

\item[(ii)]  For every $n\ge c_\theta t$,

\begin{equation}
\nonumber
\sup_{\eta\in T^{{\mathbb Z}^d},u:0\le u\le t} |S_u f(\eta)-S_u f(\eta^n)|\le
 C_\theta|||f|||_\theta e^{-\delta t},
\end{equation}
where $\eta^n$ is defined in (\ref{eta-trunc}).

\end{itemize}

\end{teo}
\begin{proof}
The proof of this theorem  is completely analogous to that of
  Theorem 3 of \cite{ramirez}. Nevertheless, in order to
mantain the content of this article self-contained, we present it here.
Note that

\begin{eqnarray}
\nonumber
&|S^n_uf(\eta)-S_uf(\eta)|=\left|\int_0^u S^n_u(L_n-L)S_{u-s}f(\eta)ds\right|\\
\nonumber
&\le D \int_0^u  \sum_{i\notin\Lambda_{n-R}}\Delta_{S_{u-s}f}
ds,
\end{eqnarray}
where $D:=\max_{i,\eta}|a_i(\eta)|+\max_{i,\eta}|b_i(\eta)|$. Therefore,
from Theorem \ref{theorem1.3}, we see that

\begin{eqnarray}
\nonumber
&
|S^n_uf(\eta)-S_uf(\eta)|\le
DA |||f|||_{\theta}\int_0^t e^{\beta (u-s)}\sum_{i\notin\Lambda_{n-R}} 
 e^{-\theta |i|_2}ds\\
\nonumber
&\le
CDA e^{(\beta -d\theta c_\theta) u}\frac{1}{1-e^{-\theta}}|||f|||_\theta, 
\end{eqnarray} 
for some constant $C$.
Choosing $c_\theta$ so that $dc_\theta -\beta>\delta$ we obtain
part $(i)$.
To prove part $(ii)$ note that for all $n\ge 1$ and $u\ge 0$
one has that for any truncation of the semigroup,
$S_u^nf(\eta)=S_u^nf(\eta^{n})$. Therefore, by the
triangle inequality we have that

$$
|S_uf(\eta)-S_uf(\eta^n)|\le
|S_uf(\eta)-S_u^nf(\eta)|+
|S_uf(\eta^n)-S_u^nf(\eta^n)|.
$$
Part $(ii)$ now follows from part $(i)$.
\end{proof}
\bigskip

\noindent Let us now continue with Gross lemma \cite{gross},
where we will need the assumption {\bf (LSI)}. It should be
noted that through the use of this assumption, here
we will be able to go from an $L_2[\nu]$ estimate to
one in the uniform norm.

\medskip

\begin{lemma}
\label{lemma3}
 {\bf [Gross lemma]}.
Consider
an infinite dimensional diffusion with  semi-group, $\{S_t:t\ge 0\}$.  Let 
$\nu$ be an invariant measure of this infinite dimensional diffusion
which satisfies {\bf (LSI)} with constant $\gamma$. Let also

\begin{equation}
\label{holder-p}
p(t)=1+e^{4t/\gamma}.
\end{equation} Then, for all $f\in L_2[\nu ]$
and $t\ge 0$, it is true that

\begin{equation}
\nonumber
||S_t f||_{p(t),\nu}\le ||f||_{2,\nu}.
\end{equation}

\end{lemma}

\begin{proof}
The proof of this lemma in the case in which the generator of
the diffusion is reversible with respect to $\nu$ can be found
in \cite{gross},\cite{holley} or in \cite{z}. To prove it in general, slight modifications
with respect to the standard proof have to be done. In particular,
it is necessary to use the invariance of $\nu$ and part $(i)$ of
Proposition \ref{prop1}.
\end{proof}

\bigskip

\noindent For each probability measure $\mu$ defined on $T^{\mathbb Z^d}$
and semi-group $\{S_t:t\ge 0\}$, we will denote by $S_t\mu$ the action
of the adjoint semi-group on $\mu$ defined by

$$
\int f dS_t\mu=\int S_t fd\mu,
$$
for every $f\in C(T^{\mathbb Z^d})$.
 Recall that for each $\Lambda\subset\mathbb Z^d$, $\nu_\Lambda$ denotes the
restriction to $T^\Lambda$ of the invariant measure $\nu$.

\medskip

\begin{lemma}
\label{lemma4}
 {\bf [Uniform norm estimate on marginal distributions]}. Consider
an infinite dimensional diffusion with truncated semi-group
at scale $n$, $\{S_t^n:t\ge 0\}$ and invariant measure $\nu$. 
Let $\eta\in T^{\Lambda_n}$, $\mu_n:=S^n_1\delta_{\eta}$ and

$$
g_{n,\eta}:=\frac{d\mu_n}{d\nu_{\Lambda_n}},
$$
  the Radon-Nikodym derivative of $\mu_n$
with respect to $\nu_{\Lambda_n}$. Then, if $q(t):=1+e^{-4t/\gamma}$, for
every $K>0$ we have that

\begin{equation}
\nonumber
\lim_{t\to\infty}\sup_\eta || g_{K t,\eta} ||_{q(s),\nu_{\Lambda_{K t}}
}\le 1,
\end{equation}
where $s$ is given by the relation $t=1+s+s^2$.
\end{lemma}

\begin{proof} Our approach will be inspired by some techniques used
in \cite{holley},
\cite{gz}
and \cite{z}. Also, we use throughout this proof the fact that
because the invariant measure $\nu$ is a Gibbs measure, its
conditional densities are positive and have bounded second order
partial derivatives.
Consider
 the semi-group  $\{S^n_t:t\ge 0\}$ of the  truncated version at scale $n$
of the infinite dimensional diffusion process with coefficients $a$ and $b$,
with infinitesimal generator $L_n$ as defined in (\ref{tg}).
Let  $P_{n,\eta}$ be the law of such a process 
starting
from $\eta\in T^{\Lambda_n}$ on the space
$C([0,\infty);T^{\Lambda_n})$ endowed with its Borel $\sigma$-algebra,
for $t\ge 0$, $P_{n,\eta,t}$ the restriction of such a law to 
 $\mathcal F_t$, the information up to time $t$ of the process.
 For each $t>0$, 
let 

$$
h_{n,\eta,t}:=\frac{d S_t^n\delta_\eta}{dm_n},
$$ 
be the Radon-Nikodym derivative of
 $S^n_t\delta_\eta$ with respect to the Lebesgue measure $m_n$ on
$T^{\Lambda_n}$.
Consider the finite dimensional diffusion on $T^{\Lambda_n}$ defined
by the infinitesimal generator

\begin{equation}
\label{efgen}
F_n:=\frac{1}{2}\sum_{i\in\Lambda_n}\partial_i ( \bar a^n_i
\partial_i ),
\end{equation}
where

\begin{equation}
\label{abarlemma}
\bar a_i^n:=E_\nu[a_i|\mathcal F_{\Lambda_n}].
\end{equation}
Let us call $Q_{n,\eta}$ the law of this diffusion starting from $\eta\in 
T^{\Lambda_n}$ defined on $C([0,\infty);T^{\Lambda_n})$ with its 
Borel $\sigma$-algebra, and for $t\ge 0$, $Q_{n,\eta,t}$ its restriction to
$\mathcal F_t$. Furthermore, call $d_{n,\eta,t}$ the Radon-Nikodym
derivative of the law of this process at time $t$ with respect to 
Lebesgue measure.
By Lemma \ref{aronson} of Appendix \ref{appendix2}, we have that



\begin{equation}
\label{ne4}
e^{-C_1|\Lambda_n|^{2}\ln|\Lambda_n|}\le d_{n,\eta,1}(\zeta)\le e^{C_1|\Lambda_n|},
\end{equation}
uniformly in $\eta,\zeta\in T^{\Lambda_n}$,
for some constant $C_1>0$.

Now, define $u_{n,1}$ as the density of
the law of the process with generator (\ref{efgen}) at time $1$
starting from $\nu_{\Lambda_n}$ with
respect to Lebesgue measure, so that

$$
u_{n,1}(\zeta):=\int d_{n,\eta,1}(\zeta)d\nu_{\Lambda_n}(\eta).
$$
It is obvious then from (\ref{ne4}) that

\begin{equation}
\label{ne4.1}
e^{-C_1|\Lambda_n|^{2}\ln|\Lambda_n|}\le u_{n,1}(\zeta)\le e^{C_1|\Lambda_n|},
\end{equation}
for all $\zeta\in T^{\Lambda_n}$.
 To prove Lemma \ref{lemma4}, let us first write

$$
||g_{K t,\eta}||^{q(s)}_{q(s),\nu_{\Lambda_{K t}}}=\int 
\left(\frac{h_{K t,\eta,1}}{d_{K t,\eta,1}}\right)^{q(s)}
\left(\frac{d_{K t,\eta,1}}{u_{K t,1}}\right)^{q(s)}
\left(\frac{u_{K t,1}}{v_{K t}}\right)^{q(s)}d\nu_{\Lambda_{K t}}.
$$
Therefore, using the bounds (\ref{ne4}) and (\ref{ne4.1}), we see that
there is a constant $C_2>0$ such that
\begin{equation}
\label{ef2}
||g_{K t,\eta}||^{q(s)}_{q(s),\nu_{\Lambda_{K t}}}\le 
\left(e^{C_2|\Lambda_{K t}|^{2}\ln |\Lambda_{K t}|}\right)^{q(s)-1}
\int 
\left(\frac{h_{K t,\eta,1}}{d_{K t,\eta,1}}\right)^{q(s)}
\left(\frac{u_{K t,1}}{v_{K t}}\right)^{q(s)-1}
\frac{d_{K t,\eta,1}}{v_{K t}} d\nu_{\Lambda_{K t}}.
\end{equation}
Let $s':=s/2$.
Now, by an application of Cauchy-Schwartz inequality to the integral
in (\ref{ef2}), we see that

\begin{eqnarray}
\nonumber
&\displaystyle
\int 
\left(\frac{h_{K t,\eta,1}}{d_{K t,\eta,1}}\right)^{q(s)}
\left(\frac{u_{K t,1}}{v_{K t}}\right)^{q(s)-1}
\frac{d_{K t,\eta,1}}{v_{K t}} d\nu_{\Lambda_{K t}}\\
\label{ef21}
&\le
\left(\displaystyle
\int 
\left(\frac{h_{K t,\eta,1}}{d_{K t,\eta,1}}\right)^{q(s)q(s')}
\frac{d_{K t,\eta,1}}{v_{K t}} d\nu_{\Lambda_{K t}}\right)^{\frac{1}{q(s')}}
\left(\displaystyle
\int 
\left(\frac{u_{K t,1}}{v_{K t}}\right)^{(q(s)-1)\frac{q(s')}{q(s')-1}}
\frac{d_{K t,\eta,1}}{v_{K t}} d\nu_{\Lambda_{K t}}
\right)^{\frac{q(s')-1}{q(s')}}.
\end{eqnarray}
Let us now see how to bound the first integral of the right-hand side of
(\ref{ef21}).
Note that 

\begin{equation}
\label{ne2}
\frac{h_{K t,\eta,1}}{d_{K t,\eta,1}}=E_{Q_{K t,\eta,1}}\left[\left. \frac{d P_{K t,\eta,1}}{d Q_{K t,\eta,1}}\right|
\mathcal F_{=1}\right],
\end{equation}
where for each $t\ge 0$, $\mathcal F_{=t}$ is the $\sigma$-algebra of events at time $t$.
Then, by Jensen's inequality and the identity (\ref{ne2}),

\begin{eqnarray}
\nonumber
&\displaystyle\int\left(\frac{h_{K t,\eta,1}}{d_{K t,\eta,1}}\right)^{q(s)q(s')}
\frac{d_{K t,\eta,1}}{v_{K t}} d\nu_{\Lambda_{K t}}
=\int\left(E_{Q_{K t,\eta,1}}\left[\left. \frac{d P_{K t,\eta,1}}{d 
Q_{K t,\eta,1}}\right|
\mathcal F_{=1}\right]\right)^{q(s)q(s')}\frac{d_{K t,\eta,1}}{v_{K t}} d\nu_{\Lambda_{K t}}\\
\label{ef1}
&\!\!\le\displaystyle
\int E_{Q_{K t,\eta,1}}\left[\left. \left(
\frac{d P_{K t,\eta,1}}{d Q_{K t,\eta,1}}\right)^{q(s)q(s')}\right|
\mathcal F_{=1}\right]\frac{d_{K t,\eta,1}}{v_{K t}} d\nu_{\Lambda_{K t}}
=E_{Q_{K t,\eta,1}}\left[ \left(
\frac{d P_{K t,\eta,1}}{d Q_{K t,\eta,1}}\right)^{q(s)q(s')}\right].
\end{eqnarray}
Now,
by the Girsanov theorem, for every natural $n$ and $t\ge 0$,

$$
\frac{d P_{n,\eta,t}}{d Q_{n,\eta,t}}=\exp\left(
\sum_{i\in\Lambda_n}\left(\int_0^t\frac{1}{\bar a^n_i}\left(\bar b^n_i-\frac{\partial \bar a^n_i}
{\partial\eta_i}\right)d\eta_i(u)-\int_0^t\frac{1}{2\bar a^n_i}\left(\bar b^n_i-
\frac{\partial \bar a^n_i}{\partial\eta_i}\right)^2 du\right)\right).
$$
where the coefficients $\bar a_i^n$ are defined in (\ref{abarlemma}) while

$$
\bar b_i^n:=E_\nu[b_i|\mathcal F_{\Lambda_n}].
$$
Therefore, again by the Girsanov theorem,  the uniform ellipticity assumption and the
boundedness of the coefficients and its derivatives, we know that there is a constant $C_3>0$ such that

\begin{equation}
\label{ef3}
E_{Q_{K t,\eta,1}}\left[ \left(
\frac{d P_{K t,\eta,1}}{d Q_{K t,\eta,1}}\right)^{q(s)q(s')}\right]
\le 
\exp\left\{C_3(q(s)^2q(s')^2-q(s)q(s'))|\Lambda_{K t}|\right\}
\le \exp\left\{2C_3e^{-4s/\gamma}|\Lambda_{K t}|\right\}.
\end{equation}
Let us now see how to bound the second term of the right-hand side
of (\ref{ef21}).
Define
$P_{n,\nu_{\Lambda_n},t}:=\int P_{n,\eta,t}d\nu_{\Lambda_n}(\eta)$
and $Q_{n,\nu_{\Lambda_n},t}:=\int Q_{n,\eta,t}d\nu_{\Lambda_n}(\eta)$
and note that $P_{n,\nu_{\Lambda_n},t}$-a.s. we have that

$$
\frac{dP_{n,\nu_{\Lambda_n},t}}{dQ_{n,\nu_{\Lambda_n},t}}(\zeta_\cdot)
=\frac{dP_{n,\zeta_0,t}}{dQ_{n,\zeta_0,t}}(\zeta_\cdot)
$$ 
Using the fact that $\nu_{\Lambda_n}$ is an
invariant measure for the diffusion with generator $L_n$,
it follows that

\begin{equation}
\label{ne211}
\frac{v_{K t}(\zeta_1)}{u_{K t,1}(\zeta_1)}= E_{Q_{K t,\nu_{\Lambda_{Kt}},1}}\left[\left.
 \frac{d P_{K t,\zeta_0,1}}{d Q_{K t,\zeta_0,1}}\right|
\mathcal F_{=1}\right].
\end{equation}
Now, for the purpose of estimating
the second term of (\ref{ef21}) we will need
to write the  Radon-Nikodym derivative in (\ref{ne211})  as

\begin{equation}
\label{ne21}
\frac{v_{K t}(\zeta)}{u_{K t,1}(\zeta)}=\int E_{Q_{K t,\eta,1}}\left[\left.
 \frac{d P_{K t,\eta,1}^{\eta'}}{d Q_{K t,\eta,1}^{\eta'}}\right|
\mathcal F_{=1}\right]d\nu_{\Lambda_{K t}}(\eta'),
\end{equation}
where $\zeta$ here plays the role of $\zeta_1$ in (\ref{ne211}),

$$
\frac{d P_{n,\eta,t}^{\eta'}}{d Q_{n,\eta,t}^{\eta'}}=\exp\left(
\sum_{i\in\Lambda_n}\left(\int_0^t\frac{1}{a^{n,\eta'}_i}\left(b^{n,\eta'}_i-\frac{\partial a^{n,\eta'}_i}
{\partial\eta_i}\right)d\eta_i(u)-\int_0^t\frac{1}{2a^{n,\eta'}_i}\left(b^{n,\eta'}_i-
\frac{\partial a^{n,\eta'}_i}{\partial\eta_i}\right)^2 du\right)\right),
$$
and

\begin{eqnarray*}
\nonumber
&a_i^{n,\eta'}(\eta(u)):=a_i^{n}(\eta(u)+\eta'-\eta)\\
&b_i^{n,\eta'}(\eta(u)):=b_i^{n}(\eta(u)+\eta'-\eta).
\end{eqnarray*}
Let

$$
r(s):=(q(s)-1)\frac{q(s')}{q(s')-1}.
$$
Then, by Jensen's inequality and the identity (\ref{ne21}), we see that

\begin{eqnarray}
\nonumber
&\displaystyle\int 
\left(\frac{u_{K t,1}}{v_{K t,1}}\right)^{r(s)}
\frac{d_{K t,\eta,1}}{v_{K t}} d\nu_{\Lambda_{K t}}(\zeta)=
\int 
\left(
\int E_{Q_{K t,\eta,1}}\left[\left. \frac{d P_{K t,\eta,1}^{\eta'}}{d 
Q_{K t,\eta,1}^{\eta'}}\right|
\mathcal F_{=1}\right]d\nu_{\Lambda_{K t}}(\eta')
\right)^{-r(s)}
\frac{d_{K t,\eta,1}}{v_{K t}} d\nu_{\Lambda_{K t}}(\zeta)\\
\nonumber
&\le
\displaystyle
\int 
\int E_{Q_{K t,\eta,1}}\left[\left. \left(
\frac{d P_{K t,\eta,1}^{\eta'}}{d 
Q_{K t,\eta,1}^{\eta'}}
\right)^{-r(s)}
\right|
\mathcal F_{=1}\right]
\frac{d_{K t,\eta,1}}{v_{K t}} d\nu_{\Lambda_{K t}}(\zeta)d\nu_{\Lambda_{K t}}(\eta')\\
\label{ef41}
&=
\displaystyle\int 
E_{Q_{K t,\eta,1}}\left[ \left(
\frac{d P_{K t,\eta,1}^{\eta'}}{d 
Q_{K t,\eta,1}^{\eta'}}
\right)^{-r(s)}
\right]
d\nu_{\Lambda_{K t}}(\eta').
\end{eqnarray}
Now, again by the Girsanov theorem,  the uniform ellipticity assumption and the
boundedness of the coefficients and its derivatives, we know that there is a constant $C_4>0$ such that

\begin{equation}
\label{ef31}
E_{Q_{K t,\eta,1}}\left[ \left(
\frac{d P_{K t,\eta,1}^{\eta'}}{d Q_{K t,\eta,1}^{\eta'}}\right)^{-r(s)}\right]
\le 
\exp\left\{C_4(r(s)-r^2(s))|\Lambda_{K t}|\right\}
\le \exp\left\{2C_4e^{-2s/\gamma}|\Lambda_{K t}|\right\}.
\end{equation}
Substituting back (\ref{ef31}) into (\ref{ef41}), 
 (\ref{ef3}) into (\ref{ef1}) and then substituting both 
(\ref{ef41}) and (\ref{ef1}) into (\ref{ef21}), we
conclude that

$$
\int 
\left(\frac{h_{K t,\eta,1}}{d_{K t,\eta,1}}\right)^{q(s)}
\left(\frac{u_{K t,1}}{v_{K t}}\right)^{q(s)-1}
\frac{d_{K t,\eta,1}}{v_{K t}} d\nu_{\Lambda_{K t}}\le
\exp\left\{C_4e^{-2s/\gamma}|\Lambda_{K t}|
\right\}.
$$
Inserting now this estimate into (\ref{ef2}) we
 see that there is a constant $C_5>0$ such that

\begin{equation}
\label{ef4}
\sup_\eta ||g_{K t,\eta}||_{q(s),\nu_{\Lambda_{K t}}}\le 
\exp\left\{C_5e^{-2s/\gamma}|\Lambda_{K t}|^{2}\ln |\Lambda_{K t}|
\right\}.
\end{equation}
Since $\lim_{t\to\infty}e^{-2s/\gamma}|\Lambda_{K t}|^{2}\ln |\Lambda_{K t}|=0$,
taking the limit when $t$ tends to $\infty$ in inequality (\ref{ef4}), we obtain Lemma \ref{lemma4}.

\end{proof}

\bigskip

\non Let us now show why do  Theorem \ref{lemma2}, Lemma \ref{lemma3}
and Lemma \ref{lemma4}, imply Theorem \ref{theorem2}. 
Let $f\in D_0$. First, note that without loss of generality, we can assume that $\langle f \rangle_\nu=0$. For $t\ge 1$, define $s\ge 0$ by the relation $t=1+s+s^2$.
Remark that by parts $(ii)$ of the truncation estimate Theorem \ref{lemma2} with $\delta=1$,  for each $\theta>0$ there exist
constants $c_\theta >0$ and $C_\theta >0$ such that,

\begin{equation}
\label{up}
|S_tf(\eta)|\le \left|S_t^{c_\theta t}f(\eta)\right|+C_\theta |||f|||_\theta e^{-t}.
\end{equation}
But, note that

$$
S_t^{c_\theta t}f(\eta)=\int S_{s+s^2}^{c_\theta t}f(\zeta^{c_\theta t})d\mu_{c_\theta t}
(\zeta^{c_\theta t}),
$$
where we use the notation $\zeta^{c_\theta t}$ for an element of $T^{\Lambda_{c_\theta t}}$,
and
where  $\mu_{c_\theta t}$ is the restriction to $\mathcal F_{\Lambda_{c_\theta t}}$ of the measure
$S^{c_\theta t}_1\delta_\eta$. Since by part $(i)$ of Theorem \ref{lemma2}
with $\delta=1$ again,
we have that uniformly in $\zeta\in T^{\mathbb Z^d}$ the
expression $S_{s+s^2}^{c_\theta t}f(\zeta)$ is exponentially close to
 $S_{s+s^2}f(\zeta)$, we conclude using (\ref{up}) that there
exist constants $c_\theta>0$ and $C_\theta >0$ such that

\begin{equation}
\label{t2}
|S_t f(\eta)|\le \left|\int S_{s+s^2} f(\zeta ) d \mu_{c_\theta t}(\zeta^{c_\theta t})\right|
+C_\theta |||f|||_\theta e^{-t},
\end{equation}
where $\zeta$ is an arbitrary extension to $T^{\mathbb Z^d}$
of $\zeta^{c_\theta t}\in T^{\Lambda_{c_\theta t}}$.
Now,


\begin{equation}
\label{t1}
\left|\int S_{s+s^2} f(\zeta)d\mu_{c_\theta t}(\zeta^{c_\theta t})\right|
=\left|\int S_{s+s^2} f(\zeta) g_{c_\theta t,\eta} d\nu_{\Lambda_{c_\theta t}}(\zeta^{c_\theta t})\right|.
\end{equation}
Let $p(t)$ be defined by (\ref{holder-p}) and $q(t)$
the conjugate exponent of $p(t)$ defined through $\frac{1}{p(t)}+\frac{1}{q(t)}=1$. Then, by H\"older's inequality, we see
that the right-hand side of inequality (\ref{t1}) is upper-bounded by

\begin{eqnarray*}
&
\displaystyle 
||S_{s+s^2} f(\zeta)||_{p(s),\nu_{\Lambda_{c_\theta t}}} 
\displaystyle ||g_{c_\theta t,\eta} 
||_{q(s),\nu_{\Lambda_{C_\theta t}}}
\le \left(||S_{s+s^2} f(\zeta)||_{p(s),\nu}
+C_\theta|||f|||_\theta e^{- t}\right) ||g_{c_\theta t,\eta} ||_{q(s),\nu_{\Lambda_{C_\theta t}}}\\
&
\le\left(||S_{s^2} f(\zeta)||_{2,\nu}+C_\theta |||f|||_\theta e^{- t}\right) 
||g_{c_\theta t,\eta} ||_{q(s),\nu_{\Lambda_{C_\theta t}}}\\
&\le
\left(e^{-as^2/\gamma} ||f||_{2,\nu}+C_\theta |||f|||_\theta e^{- t}\right) 
\sup_\eta||g_{c_\theta t,\eta} ||_{q(s),\nu_{\Lambda_{C_\theta t}}},
\end{eqnarray*}
where in the first inequality we have used the truncation estimate of Theorem \ref{lemma2},
in the second inequality Gross lemma (Lemma \ref{lemma3}),
 and in the last
inequality  Corollary \ref{lemma1} with spectral gap constant $g=1/\gamma$.
Now, by Lemma \ref{lemma4},
there is a $t_0>0$ such that 

$$K_1(\nu,\theta):=\sup_{t\ge t_0}\sup_\eta ||g_{c_\theta t,\eta} ||_{q(s),\nu_{\Lambda_{C_\theta t}}}
<\infty.
$$ Thus, by inequality (\ref{t2}), and the previous development,
for $t\ge t_0$ we have that

\begin{equation}
\label{endend}
\left|S_t f(\zeta)\right|\le K_2(\nu,\theta)\left(e^{-a s^2/\gamma} ||f||_{2,\nu}+
C_\theta |||f|||_\theta e^{- t}\right),
\end{equation}
where $K_2(\nu,\theta):=K_1(\nu,\theta)+1$. Now, note that

\begin{equation}
\label{l2triple}
||f||_{2,\nu}\le \sqrt{\gamma\sum_{i\in\mathbb Z^d}\sup_{\eta}\left(
\partial_i f\right)^2}\le 
\sqrt{\gamma\sum_{i\in\mathbb Z^d}\sup_{\eta}\left(
\partial_i^2 f\right)^2}\le 
\gamma^{1/2}\sqrt{|||f|||_\theta},
\end{equation}
where in the first inequality we have used the spectral gap inequality
(\ref{thespectral}),
while in the second, the fact that $\sup_{\eta}(\partial_i f)
\le \sup_{\eta}(\partial^2_i f)$ (which is a consequence of
the identity $\partial_i f(\eta)=\int_{\eta'_i}^{\eta_i}\partial_i^2
f(\zeta)d\zeta_i$, where $\zeta_j=\eta_j$ for $j\ne i$, valid for all $\eta\in T^{\mathbb Z^d}$ and a
$\eta'_i$ depending on $\eta$). Substituting (\ref{l2triple}) into
(\ref{endend}), we conclude the proof of Theorem \ref{theorem2}.

\smallskip

\section{Proof of Theorem \ref{theorem1}.} 

 We need to recall some basic notions
(see for example Kato \cite{kato}) which will be used throughout this section.
Given a closed operator $T$ defined on $L_2[\nu]$, we will
denote its spectrum by $\Sigma(T)$. We will
denote by $R(z,T)$ the resolvent operator for every $z\notin\Sigma(T)$.
 We wil say that the spectrum $\Sigma(T)$ contains
a bounded part $\Sigma_2\ne \emptyset$ separated from the rest $\Sigma_1\ne\emptyset$ if
there exists a rectifiable, simple closed curve $\Gamma$ 
in the complex plane $\mathbb C$ which encloses an open set
containing $\Sigma_2$ in its interior and
$\Sigma_2$ in its exterior. Of course
we have $\Sigma(T)=\Sigma_1\cup\Sigma_2$.
Throughout, given a vector space $X$ and two 
subspaces $Y$ and $Z$ we will use the standar notation
for direct sum $X=Y\oplus Z$, whenever for every
vector $x\in X$ we have a unique decomposition
$x=y+z$ for some $y\in Y$ and $z\in Z$.
We will say that two subspaces $M_1$ and $M_2$ of $L_2[\nu]$
form a {\it decomposition associated to $\Sigma_1$ and $\Sigma_2$} if
$L_2[\nu]=M_1\oplus M_2$, the spectrum of $T_{M_1}$ is $\Sigma_1\cup\{0\}$ and the
spectrum of $T_{M_2}$ is $\Sigma_2\cup\{ 0\}$. Here $T_{M_1}:=P_1T=TP_1$ and $T_{M_2}:=P_2T=TP_2$,
where 

$$
P_1=\frac{1}{2\pi i}\int_\Gamma  R(z,T)dz,
$$
 is the projection of $L_2[\nu]$ onto $M_1$ along $M_2$ and

\begin{equation}
\label{proj}
P_2={\mathbb I}-P_1
\end{equation}
the projection of $L_2[\nu]$ onto $M_2$ along $M_2$.

Firstly, we derive the following
 proposition giving some basic information about 
the character of $0$ in the spectrum of the unperturbed operator.

\medskip
\begin{proposition}
\label{prop2} 
Consider
an infinite dimensional  diffusion with coefficients $a$ and $b$
and infinitesimal generator $(L_0,D(a,b))$.
Let $\nu$ be an invariant Gibbs  measure of this infinite dimensional
diffusion satisfying {\bf (LSI)} and $(L_0,\bar D(a,b))$ be its infinitesimal
generator on $L_2[\nu]$. 
Then the following statements are satisfied.

\begin{itemize}
\item[(i)]  $0$ is a  simple eigenvalue
of the operator $(L_0, \bar D(a,b))$.

\item[(ii)] The intersection of the
open disc centered at $0$ of radius $a/(4\gamma)$ with the spectrum
$\Sigma(L_0)$ of $L_0$ is $\{0\}$, so that $0$ is an isolated eigenvalue.
\end{itemize}

\end{proposition}
\begin{proof} {\it Part $(i)$.} 
 It is enough to prove
that $0$ is a simple eigenvalue of the adjoint
operator $(L_0^*,\bar D(a,b)^*)$. Assume that $g\in \bar D(a,b)$ is a
normalized function such that

\begin{equation}
\label{rs3}
\int (L_0^*g)fd\nu=0,
\end{equation}
for every $f\in \bar D(a,b)$. To prove that $0$ is a simple eigenvalue
of $L_0^*$, it is enough to show that this implies that
$\nu$-a.s. $g=1$. But the left-hand side of (\ref{rs3}) can
be written as $\int g L_0 f d\nu=0$. Since $\nu$ is
a Gibbs measure that satisfies
the logarithmic Sobolev inequality, by Corollary \ref{cor}, $\nu$ is the
unique invariant measure. It follows that necessarily $\nu$-a.s. $g=1$.

\smallskip

\noindent {\it Part $(ii)$.}
From part $(iii)$ of
Proposition \ref{prop1}, we see that for every complex $z$ such that 
$0<|z|<\frac{a}{2\gamma}-|z|$
one has that for every $f\in \bar D(a,b)$,

$$
\| (L_0-z)f\|_{2,\nu}\ge m \|f\|_{2,\nu},
$$
where $m:=\min\{|z|,\frac{a}{2\gamma}-|z|\}$. This shows that 
every $z$ such that $0<|z|<\frac{a}{4\gamma}$ is in the resolvent set of $L_0$.

\end{proof}

\smallskip

\noindent Let us now recall that for each real $\epsilon$, 

\begin{equation}
\label{eleps}
L_\epsilon=L_0+\epsilon A,
\end{equation}
where $A=\sum c_i\partial_i$ is a diagonal first order operator
satisfying condition {\bf (R)} with respect to the invariant measure $\nu$.
Note that by  Corollary \ref{coro0} we know that
 for every real $\epsilon$, 
the operator $(L_\epsilon,D_0)$
is closable, having the same domain $\bar D(a,b)$ as the closure
 of $(L_0,D_0)$ on $L_2[\nu]$.

Let $\epsilon_0>0$.
We  say that a family $\{T(\epsilon):\epsilon\in (-\epsilon_0,\epsilon_0)\}$ 
   of bounded operators defined on $L_2[\nu]$, is {\it holomorphic}
in $\epsilon$ if and only if each $\epsilon$ has a neighborhood
in which $T(\epsilon)$ is bounded and $(f,T(\epsilon)g)_\nu$ is
holomorphic for every $f,g$ in a dense subset of $L_2[\nu]$ (see section VII.1
of Kato \cite{kato}). Recall the definition of $\epsilon_c$
in (\ref{epsc}).

\smallskip

\begin{lemma}
\label{lemmaf} 
Consider
an infinite dimensional diffusion with coefficients $a$ and $b$
and infinitesimal generator $(L_0,D(a,b))$.
Let $\nu$ be an invariant Gibbs measure of this infinite dimensional
diffusion satisfying {\bf (LSI)} and $(L_0,\bar D(a,b))$ be its infinitesimal
generator on $L_2[\nu]$. 
Let $L_\epsilon$ be given by (\ref{eleps}), with $A$ a diagonal
first order operator satisfying condition {\bf (R)} with respect to $\nu$.
 Consider the complex contour $\Gamma:=\{z\in{\mathbb C}:
|z|=a/(4\gamma)\}$. 
Then, the projection $P_{\epsilon}$ of $L_2[\nu]$
onto $M_{1,\epsilon}$ along $M_{2,\epsilon}$ 
is holomorphic 
as a function of $\epsilon$ for $\epsilon\in(- \epsilon_c,\epsilon_c)$.
 and admits the following expansion with radius of convergence $\epsilon_c$

$$
P_\epsilon=-\frac{1}{2\pi i}\sum_{k=0}^\infty\epsilon^k
\int_\Gamma R(z,L_0) 
(- A\ R(z,L_0))^kdz.
$$
\end{lemma}
\begin{proof}
Note that for $z\in\Gamma$, whenever $\epsilon$ 
is non-negative and is such that
$\|\epsilon A\ R(z,L_0)\|_{2,\nu}<1$, we have the following
expansion (see Theorem 1.5, page 66 and chapter VIII of \cite{kato})

$$
R(z,L_\epsilon)=R(z,L_0)\sum_{k=0}^\infty \epsilon^k(- A\ R(z,L_0))^k.
$$
Observe that for $z\in \Gamma$, $\left\|\frac{L_0}{L_0-z}\right\|_{2,\nu}\le 2$
$\| R(z,L_0)\|_{2,\nu}\le \frac{\gamma}{2a}$.
Therefore, we can apply Corollary \ref{coro0} to
conclude that for every $\lambda>0$,

$$
\|\epsilon A\ R(z,L_0)\|_{2,\nu}\le\frac{2\epsilon C_0}{\sqrt{a}}
\frac{1}{\lambda}
+\frac{\epsilon C_0\gamma}{2a^{3/2}}\lambda.
$$
Taking the infimum over $\lambda>0$, we conclude that

$$
\|\epsilon A\ R(z,L_0)\|_{2,\nu}\le\epsilon\frac{C_0\sqrt{\gamma}}{a},
$$
which proves the analyticity of the resolvent operator $R(z,L_\epsilon)$
for $z\in\Gamma$ when $\epsilon\in (-\epsilon_c,\epsilon_c)$.
Finally,  by (\ref{proj}), the projection $P_{\epsilon}$ of $L_2[\nu]$ onto
$M_{1,\epsilon}$ along $M_{2,\epsilon}$ can be expressed as

$$
P_{\epsilon}=-\frac{1}{2\pi i}\int_\Gamma R(z,L_\epsilon)dz,
$$
which proves the lemma.

\end{proof}

\smallskip
From Lemma \ref{lemmaf}, we have directly the following corollary regarding
the adjoints $L_0^*$ and $A^*$ of $L_0$ and $A$ in $L_2[\nu]$
respectively.
\smallskip

\begin{coro}
\label{coro4}
 Let
$A,L_0,L_\epsilon$ and $\nu$ be as in Proposition
\ref{prop2}. Consider the complex contour $\Gamma:=\{z\in{\mathbb C}:
|z|=a/4(\gamma)\}$. 
Then, 
the projection $P^*_{\epsilon}$ of $L_2[\nu]$
onto $M^*_{1,\epsilon}$ along $M^*_{2,\epsilon}$ 
is holomorphic 
as a function of $\epsilon$ for $\epsilon\in(- \epsilon_c,\epsilon_c)$.
 and admits the following expansion with radius of convergence $\epsilon_c$

$$
P^*_\epsilon=-\frac{1}{2\pi i}\sum_{k=0}^\infty\epsilon^k
\int_{\bar\Gamma} 
(-  R(z,L^*_0)A^*)^kR(z,L^*_0) dz.
$$
\end{coro}
\smallskip

\non Let us now prove Theorem \ref{theorem1}. By part $(i)$  of 
Proposition \ref{prop2} and Corollary \ref{coro4}, we 
see that for each $\epsilon\in (-\epsilon_c,\epsilon_c)$ there exists
a unique invariant measure $\nu_\epsilon$ of the
infinite dimensional diffusion with generator $L_\epsilon$ 
in $M_2[\nu]$. On the other hand, we know that $g:=1$ is
an eigenfunction associated to the eigenvalue $0$  of $L_0$ in  $L_2[\nu]$. Let

$$
g'_\epsilon:=P^*_\epsilon g.
$$
By part $(iii)$ of Corollary \ref{coro4}, we know that $g'_\epsilon$
admits the expansion

$$
g'_\epsilon=\sum_{k=0}^\infty \epsilon^k f'_k,
$$
where $f'_0:=g$ and

$$
f'_k:=-\frac{1}{2\pi i}
\int_{\bar\Gamma} 
(-  R(z,L^*_0)A^*)^kR(z,L^*_0)g dz.
$$
By parts $(i)$ and $(ii)$ of Corollary \ref{coro4}, necessarily
$L_\epsilon^*g'_\epsilon=0$. Hence, for every $f\in \bar D(a,b)$,

$$
\sum_{k=0}^\infty \epsilon^k (f'_k,(L_0+\epsilon A)f)_\nu=0.
$$
Matching equal powers of $\epsilon$ in the above equation, we
conclude that for each $k\ge 0$, $h:=f_{k+1}$ is solution of
the equation

\begin{equation}
\label{ss}
L_0^*h=-A^*f'_k.
\end{equation}
Since the kernel $ker(L^*_0)$ of the operator $L_0^*$ is one-dimensional and
$A^*f'_k$ is orthogonal to $ker(L^*_0)$, it follows that 
 the sequence of functions $f_0:=f'_0$ and
$f_k:=f'_k-\langle f'_k\rangle_\nu$, $k\ge 1$,
is the only sequence satisfying (\ref{ss})
under the condition that the average of each term with respect to $\nu$
vanishes. But since $\langle g_\epsilon\rangle_\nu=0$, we
see that $\left\langle\sum_{k=1}^\infty\epsilon^kf'_k\right\rangle_\nu=0$.
It follows that
$g_\epsilon:=1+\sum_{k=1}^\infty \epsilon^k f_k$
is the Radon-Nikodym derivative of $\nu_\epsilon$ with
respect to $\nu$.
  This finishes the
proof of Theorem \ref{theorem1}. 
\smallskip

\section{Interacting diffusions.} 

 Here we will consider
the effect of an interaction given by a diagonal first order operator
satisfying condition {\bf (R)}
 on a set of independent diffusions.
Consider a set of coefficients $a,b:T^{\mathbb Z^d}\to [0,\infty)$,
with $a=\{a_i:i\in\mathbb Z^d\}$ and $b=\{b_i:i\in\mathbb Z^d\}$
such that $a$ is uniformly elliptic. We will also assume that
for each $i\in\mathbb Z^d$, $a_i$ and $b_i$ are functions only
of $\eta_i$. Obviously, the infinitesimal generator $L_0$
with these coefficients defines the dynamics of independent
diffusion processes indexed by $\mathbb Z^d$.
We say that $L_0$ is the infinitesimal generator of
a set of {\it independent diffusions}. Therefore
the invariant measure of the corresponding infinite dimensional
diffusion is a product measure $\nu$. Obviously $\nu$ satisfies
the logarithmic Sobolev inequality. We want to quantify the
effect of a perturbation on the dynamics which introduces interaction
on the invariant measure. Here we will prove the following result.

\medskip

\begin{teo}
\label{theorem3} Consider the generator $L_0$ of a set of independent diffusions
with invariant product measure $\nu$.
Let $A$ be a diagonal first order operator satisfying condition {\bf (R)} with respect to $\nu$, with coefficients which are smooth.  Then, for each
$\epsilon\in (-\epsilon_c,\epsilon_c)$, the diffusion with infinitesimal
generator

$$
L_\epsilon=L_0+\epsilon A,
$$
has a unique invariant measure $\nu_\epsilon$ which has a Radon-Nikodym
derivative $g_\epsilon$ with respect to $\nu$ of the form

\begin{equation}
\label{expoh}
g_\epsilon=1+\sum_{k=1}^\infty\epsilon^k f_k,
\end{equation}
where for each $k\ge 1$,

$$
f_k=\sum_{i_1,\ldots,i_k\in\mathbb Z^d}f_{i_1,\ldots,i_k}^{(k)},
$$
and for each $i_1,\ldots,i_k\in\mathbb Z^d$,
$f_{i_1,\ldots,i_k}^{(k)}\in L_2[\nu]$ is a local function depending only
on $\cup_{j=1}^k\Lambda(i_j,R)$. Furthermore, there
is a constant $C>0$ such that for each $i,j\in\mathbb Z^d$, we have that

\begin{equation}
\label{nunu}
\sup_{A,B\in\mathcal B(T)}\left|\nu_\epsilon(\eta_i\in A,\eta_j\in B)-\nu_\epsilon(\eta_i\in A)\nu_\epsilon(\eta_j\in B)\right|
\le  
\begin{cases}
\epsilon\qquad &{\rm if}\ 1\le |i-j|_2\le R\\
\epsilon^2\qquad &{\rm if}\ |i-j|_2> R
\end{cases}.
\end{equation}

\end{teo}

\smallskip

It is important to remark that the terms in the expansion
(\ref{expoh}) of Theorem \ref{theorem3}, can be in principle
computed explicitly. For example, in the case in which
the independent diffusions are Brownian motions with drift,
so that at site $i\in\mathbb Z^d$, we have a Brownian
motion with drift $v_i$,
we can easily check that

$$
f_i^{(1)}=\int_0^\infty \int_{T^{\Lambda(i,R)}}p(t,\eta_{\Lambda(i,R)},\zeta)
\partial_i c_i(\zeta)dm_{\Lambda(i,R)}(\zeta)dt,
$$
where

$$
p(t,\zeta',\zeta):=\prod_{j\in\Lambda(i,R)}p_j(t,\zeta'_j,\zeta_j),
$$
and $p_j(t,\zeta'_j,\zeta_j)$ is the transition probability
density of a  Brownian motion with drift $-v_j$ on $T$.

\medskip

Let us now prove Theorem \ref{theorem3}.
By Theorem \ref{theorem1}, for $\epsilon<\epsilon_c$, 
we know that the diffusion with generator 
$L_\epsilon:=L_0+\epsilon A$, has an invariant measure
with a Radon-Nikodym derivative $g_\epsilon$ with respect
to $\nu$, which admits the expansion (\ref{expoh}).
Furthermore, Theorem \ref{theorem1}, also tells us that
 $L^*_0f_{k+1}=-A^*f_k$, $\langle f_k\rangle_\nu=0$ for $k\ge 0$ and $f_0:=1$. In particular, $f_1$
satisfies the equation

$$
L_0^*f_1=\sum_i 
\partial_i (c_i v_i),
$$
where $v_i$ is the Radon-Nikodym derivative with respect
to the Lebesgue measure of the
marginal at site $i$ of the invariant measure $\nu$. From here we conclude that

\begin{equation}
\nonumber
f_1(\eta)=\sum_i f_i^{(1)}(\eta),
\end{equation}
where  $f_i^{(1)}$ is the solution of the equations

\begin{eqnarray}
\label{ef11}
&L_0^*f_i^{(1)}(\eta)=\partial_i (c_i v_i)\\
\nonumber
&\langle f_i^{(1)}\rangle_\nu=0.
\end{eqnarray}
From (\ref{ef11}) we can see that for each $i\in\mathbb Z^d$,
$f_i^{(1)}$ is a local function depending only on the sites on
 $\Lambda(i,R):=\{j\in\mathbb Z^d:|i-j|_2\le R\}$, so that
$f_1$ is a sum of local functions of this type. 
For the terms of order $\epsilon^2$, we can see that

\begin{equation}
\nonumber
f_2(\eta)=\sum_{i,j}f_{i,j}^{(2)}(\eta),
\end{equation}
where this time, for each $i,j\in\mathbb Z^d$, $f_{i,j}^{(2)}$ is the
solution of

\begin{eqnarray}
\nonumber
&L_0^*f_{i,j}^{(2)}(\eta)=\partial_i(c_iv_i f_j^{(1)})\\
\nonumber
&\langle f_{i,j}^{(2)}\rangle_\nu=0.
\end{eqnarray}
From here we see that for each $i,j\in\mathbb Z^d$,
 $f_{i,j}^{(2)}$ is a local function depending only on the
coordinates of $\Lambda(i,R)\cup\Lambda(j,R)$.
By a similar argument we can conclude that for each $k\ge 1$,

$$
f_k=\sum_{i_1,\ldots,i_k\in\mathbb Z^d}f_{i_1,\ldots,i_k}^{(k)},
$$
where for each $i_1,\ldots,i_k\in\mathbb Z^d$, $f_{i_1,\ldots,i_k}^{(k)}$
is a local function depending only on the coordinates of $\cup_{l=1}^k
\Lambda(i_l,R)$.
From the above properties of the
expansion (\ref{expoh}), we can deduce (\ref{nunu}).

\medskip

\appendix
\section{Proof of Theorem \ref{theorem1.3}}
\label{appendix}
The proof of Theorem \ref{theorem1.3} is a slight modification
of Theorem 2 of \cite{ramirez}. Nevertheless, in order to give 
 self-contained
proofs of Theorems \ref{theorem2} and \ref{theorem1},
 have decided to include
the details here. We recall that $R$ denotes the range of the coefficients
$a$ and $b$ of ther diffusion.

Consider the space of functions

$$
D(T^{\mathbb Z^d}):=\{g\in C^2(T^{\mathbb Z^d}):|g|_\Delta:=
\sum_{i\in\mathbb Z^d}\Delta_g(i)<\infty\}.
$$
For each $n\ge 1$ define
the coefficients $a^n$ and $b^n$, by
$a^n_i:=(a_i-1)\mathbbm 1_{\Lambda_n}(i)+1$ and
$b^n_i:=b_i \mathbbm 1_{\Lambda_n}(i)$,  and the generator
$(\bar L_n, D_n(a,b))$ with coefficients $a^n$ and $b^n$. For $f\in C^2(T^{\mathbb Z^d})$
we consider the solution $u^n(t,\eta)$ of the
equation

\begin{eqnarray}
\nonumber
& \frac{\partial u^n}{\partial t}=\bar L_n u^n \quad {\rm for}\ t>0,\\
\nonumber
& u^n(0,\eta)=f(\eta).
\end{eqnarray}
Note that since the coefficients of the operator $\bar L_n$
are differentiable with derivatives which are H\"older continuous of
positive parameter, 
by classical regularity theory for parabolic artial differential
equations (see for example \cite{varadhan})
we can define for $i,j\in\mathbb Z^d$, $u^n_i:=\frac{\partial u^n}{\partial\eta_i}$
and $u^n_{i,j}:=\frac{\partial^2 u^n}{\partial\eta_i\partial\eta_j}$.
We will show that for each $\vartheta>0$ there are constants $A$
and $\beta$ (not depending on $n$ nor $f$) such that

\begin{equation}
\label{uenei}
|u^n_{i,j}(t,\eta)|\le A|||f|||_\vartheta e^{\beta t}e^{-\frac{\vartheta}{2}(|i|_2+|j|_2)}.
\end{equation}
Note that $u^n_j$ is solution of the equations

\begin{eqnarray}
\nonumber
& \frac{\partial u^n_j}{\partial t}=\bar L_n u^n_j
+\sum_{i\in\mathbb Z^d}\left(\frac{1}{2}\frac{\partial \bar a_i}{\partial\eta_j}
\frac{\partial u^n_i}{\partial\eta_i}
+\frac{\partial \bar b_i}{\partial\eta_j}u_i^n\right)
 \quad {\rm for}\ t>0,\\
\nonumber
& u^n_j(0,\eta)=\frac{\partial f}{\partial\eta_j}(\eta).
\end{eqnarray}
Now define for $t\ge 0$ and $\eta\in T^{\mathbb Z^d}$,

$$
W^n(t,\eta):=\sum_{j\in\Lambda_n} e^{2\vartheta |j|_2} (u_j^n)^2(t,\eta).
$$
By arguments similar to the proof of Theorem 2 of \cite{ramirez},
we conclude that

\begin{eqnarray}
\nonumber
&\frac{\partial W^n}{\partial t}-\bar L_n W^n\le
K_1 W^n\quad  {\rm for}\ t>0\ {\rm and}\ \eta\in T^{\mathbb Z^d}\\
\nonumber
&W^n(0,\eta)=\sum_{i\in\Lambda_n}e^{2\vartheta |i|_2}\left(\frac{\partial f(\eta)}{
\partial\eta_i}\right)^2,
\end{eqnarray}
where

$$
K_1:=\left(\frac{(A')^2}{4a}(2R+1)^d+B'\right)(2R+1)^d e^{2\vartheta R}
+(2R+1)^d B',
$$
and where we recall the definition of $a$ given in (\ref{uel}),

\begin{eqnarray*}
&A':=\sup_{i,j\in\mathbb Z^d,\eta\in T^{\mathbb Z^d}}\left|\frac{\partial a_i}{\partial
\eta_j}\right|\quad {\rm and}\\
&B':=\sup_{i,j\in\mathbb Z^d,\eta\in T^{\mathbb Z^d}}\left|\frac{\partial b_i}{\partial
\eta_j}\right|.
\end{eqnarray*}
From here, using the fact that since $|\partial_i f|\le |\partial_i^2 f|$,
one has that$|W^n(0,\eta)|\le |||f|||_{\vartheta}$, it follows that

\begin{equation}
\label{win}
W^n(t,\eta)\le e^{K_1 t}
|||f|||^2_{\vartheta}.
\end{equation}
Next define

$$
V^n(t,\eta):=\sum_{j,k\in\Lambda_n}(u^n_{j,k})^2 e^{\vartheta (|k|_2+|j|_2)}.
$$
As in Theorem 2 of \cite{ramirez} we conclude this time that

\begin{eqnarray}
\nonumber
&\frac{\partial V^n}{\partial t}-\bar L_n V^n\le
K_2 V^n+K_3 W^n\quad  {\rm for}\ t>0\ {\rm and}\ \eta\in T^{\mathbb Z^d}\\
\label{vvin}
&V^n(0,\eta)=\sum_{j,k\in\Lambda_n}e^{\vartheta (|k|_2+|j|_2)}
\left(\frac{\partial^2 f(\eta)}{
\partial\eta_k\partial\eta_j}\right)^2,
\end{eqnarray}
where

\begin{eqnarray}
\nonumber
& K_2:=e^{\vartheta R}(2R+1)^{2d}\left(
\frac{(A')^2}{4a}+2B'+
\frac{A''}{2}e^{\vartheta R}
\right)
+(2R+1)^d \left(2B'+\frac{A''}{2}+B''\right)\quad{\rm and}\\
\nonumber
& K_3:=B'' e^{2\vartheta R}(2R+1)^{2d},
\end{eqnarray}
with

\begin{eqnarray*}
&A'':=\sup_{i,j,k\in\mathbb Z^d,\eta\in T^{\mathbb Z^d}}\left|\frac{\partial^2 a_i}{\partial
\eta_j\partial\eta_k}\right|\quad {\rm and}\\
&B'':=\sup_{i,j,k\in\mathbb Z^d,\eta\in T^{\mathbb Z^d}}\left|\frac{\partial^2 b_i}{\partial
\eta_j\partial\eta_k}\right|.
\end{eqnarray*}
It follows now from (\ref{win}) and (\ref{vvin}) that

\begin{equation}
\nonumber
V^n(t,\eta)\le e^{K_2 t} |||f|||^2_{\vartheta}
+\frac{K_3}{K_1+K_2}e^{2(K_1+K_2)t}|||f|||^2_{\vartheta}.
\end{equation}
It follows that (\ref{uenei}) is satisfied with

$$
A^2:=\left(\frac{K_1+K_2+K_3}{K_1+K_2}\right).
$$
and

$$
\beta:=K_1+K_2.
$$
Following the argument given in the proof of Theorem 3 of \cite{ramirez}
using the theorem of Arzel\'a-Ascoli and the fact that any uniformly
convergent subsequence of $\{u^n:n\ge 0\}$ to a function $u$ is a solution
of (see \cite{stroock})

\begin{eqnarray}
\nonumber
& \frac{\partial u}{\partial t}=\bar L_n u^n \quad {\rm for}\ t>0,\\
\nonumber
& u(0,\eta)=f(\eta),
\end{eqnarray}
so that for all $i,j\in\mathbb Z^d$ it must be true that

\begin{equation}
\nonumber
|u_{i,j}(t,\eta)|\le A|||f|||_\vartheta e^{\beta t}e^{-\frac{\vartheta}{2}(|i|_2+|j|_2)},
\end{equation}
we conclude the proof of the theorem.
\bigskip

\section{Heat kernel bounds}
\label{appendix2}
Here we will prove the following lemma about heat kernel bounds.

\medskip

\begin{lemma}
\label{aronson} Let $N\ge 1$ and $a_1,\ldots, a_N$ a set of
real functions  defined on $T^N$ with bounded second order
partial derivatives. Consider a finite dimensional diffusion with
generator

$$
L=\frac{1}{2}\sum_{j=1}^N\frac{\partial}{\partial\eta_j}\left(
a_j (\eta)\frac{\partial}{\partial\eta_j}\right).
$$
Let $\Gamma_{N,\eta,t}$ be the Radon-Nikodym derivative of the law
of this diffusion starting from $\eta\in T^N$ at time $t$ with respect to the Lebesgue measure. Then, there exists a constant $C_1$ such that
the following estimates hold.

$$
e^{-C_1 N^{2}\ln N}\le \Gamma_{N,\eta,1}(\zeta)\le e^{C_1 N}.
$$

\end{lemma}
\begin{proof} Let us first prove the upper bound. We will
essentially use ideas of Nash \cite{n} as presented in
Fabes-Stroock \cite{fs}. Lemma 1 of \cite{ramirez}
(see also Proposition 1 of page 30 of \cite{lu}) gives
us the bound for all $\eta,\zeta\in T^N$

\begin{equation}
\label{fs1}
\Gamma_{N,\eta,1}(\zeta)\le (u_2)^N\Gamma_0(u_1,\zeta,\eta),
\end{equation}
for constants $u_1$ and $u_2$ which do not depend on $N$,
where $\Gamma_0$ is the fundamental solution of $\frac{1}{2}\Delta-
\frac{\partial}{\partial t}$ on $T^N$.
On the other hand, for every $t>0$ and $\eta,\zeta\in T^N$ we have that

$$
\Gamma_0(t,\zeta,\eta)=\sum_{j\in\mathbb Z^N}
e^{-2\pi^2|j|_2t} e^{2\pi i j\cdot (\eta-\zeta)}
\le \sum_{j\in\mathbb Z^N}
e^{-2\pi^2|j|_2t}\le\left(1+\sqrt{\frac{2}{t}}\right)^N.
$$
Substituting this back into (\ref{fs1}) with $t=1$, we obtain the upper bound.
To prove the lower bound,
 denote by $\pi: T^N\to\mathbb [0,1]^N$ the
bijection which assigns to each $\eta\in T^N$ its representative
in $[0,1]^N$.
Now note that

$$
\Gamma_{N,\eta,t}(\zeta)=\sum_{i\in\mathbb Z^d} \Gamma'_{N,\eta,t}(\pi(\zeta)+i),
$$
where $\Gamma'_{N,\eta,t}$ is the fundamental solution of a diffusion
in $\mathbb R^N$ with infinitesimal generator 

$$
L'=\frac{1}{2}\sum_{j=1}^N\frac{\partial}{\partial x_j}\left(
a'_j(x)\frac{\partial}{\partial x_j}\right).
$$
where  $a'_j(x):=a_j(\pi^{-1}(x))$. Therefore to prove the lower bound
it is enough to show that for every $x\in [0,1]^N$
one has that

$$
\Gamma'_{N,0,1}(x)\ge e^{-C_1 N^{2}\ln N}.
$$
From Lemma 2.6 of \cite{fs}, we know that in the case $|x|_2\le 1$,
one has that

$$
 \Gamma'_{N,0,1}(x)\ge \frac{1}{C_1'},
$$
for some constant $C_1'$ depending only on $a$ [c.f. (\ref{uel})].
Therefore, we can assume without loss of generality that
$|x|_2> 1$. We will now follow the argument presented in the
proof of Theorem 2.7 of \cite{fs}. Let $k$ be the smallest integer
larger than $4|x|_2$. Let $\chi_1,\ldots,\chi_{k-1}$ be
such that $|jx/k-\chi_j|_2< 1/(2\sqrt{k})$ for $1\le j\le k-1$,
$|\chi_1|_2<1/\sqrt{k}$ and $\max_{1<j<k}|\chi_j-\chi_{j-1}|_2<1/\sqrt{k}$.
Then, by Lemma 2.6 of \cite{fs} we conclude that there is a
constant $C_1''$ such that

\begin{eqnarray*}
&\Gamma'_{N,0,1}(x)=\int \Gamma'_{N,1/k,0}(\chi_{1}) \Gamma'_{N,1/k,\chi_{1}}(\chi_{2})
\cdots \Gamma'_{N,1/k,\chi_{k-1}}(x)d\chi_1\cdots\chi_{k-1}\\
&\ge \frac{C_1''}{k^{N/2}}\left(\frac{B_{N}}{2^{N/2}C_1''}\right)^k\ge
\frac{C_1''}{|x|_2^{N}}\left(\frac{B_{N}}{2^{N/2}C_1''}\right)^{4|x|_2^2}.
\end{eqnarray*}
where $B_N$ is the Lebesgue measure of an $N$-dimensional
Euclidean ball of radius $1$. But $|B_N|\ge\left(\frac{1}{N}\right)^{N/2}$.
On the other hand since $x\in [0,1]^N$, we have that $|x|_2\le \sqrt{N}$.
Therefore

$$
\Gamma'_{N,0,1}(x)\ge \frac{C_1''}{N^{N/2}}
\left(\frac{1}{(2N)^{N/2}C_1''}\right)^{4N}\ge e^{-C_1 N^2\ln N},
$$
which proves the lemma.

\end{proof}

\medskip
\noindent {\bf Acknowledgments}: The author thanks Olivier Bourget
for useful discussions related to Theorem \ref{theorem1}
and the referee for detailed comments on a previous version of
this paper.

\nopagebreak


\begin{thebibliography}{99}
\addcontentsline{toc}{chapter}{Bibliography}

\bibitem{cr15} D. Campos and A.F. Ram\'\i rez, {\em Asymptotic expansion of the invariant measure for ballistic random walk in the low disorder regime},  arXiv:1511.02945 (2015).

\bibitem{dynkin} E. Dynkin, {\em Markov Processes}, I. Academic Press, New York (1965).

\bibitem{fs} E. Fabes and D. Stroock, {\em A new proof of Moser's
   parabolic Harnack inequality using the old ideas of Nash}, J. Functional
  Analysis {\bf 42}, 29-63 (1981).


\bibitem{gt} D. Gilbarg and N.S. Trudinger, {\em Elliptic partial
differential equations of second order}, Grundlehren der Mathematischen Wissenschaften 224, Springer-Verlag, Berlin (1983).

\bibitem{gross} L. Gross, {\em Logarithmic Sobolev inequalities}, 
Amer. J. Math. {\bf 97}, 1061-1083 (1975).

\bibitem{gz} A. Guionnet and B. Zegarlinski,
 {\em Lectures on logarithmic Sobolev inequalities},  S\'eminaire de
 Probabilit\'es,  XXXVI,  1-134, Lecture Notes in Math.,
{\bf  1801}, Springer, Berlin, (2003).


\bibitem{holley} R. Holley, {\em The one-dimensional stochastic $X-Y$ model},
Collection: Random walks, Brownian motion, and interacting particle systems,
Progr. Probab. {\bf 28}, Birkhauser, 295-307 (1991).


\bibitem{stroock} R. Holley and D. Stroock, {\em Diffusions on an Infinite
Dimensional Torus}, J. Functional Analysis {\bf 42},29-63 (1981).

\bibitem{kato} T. Kato, {\em Perturbation theory for linear operators},
Springer-Verlag, Berlin Heidelberg (1995)


\bibitem{ko} T. Komorowski and S. Olla, {\em 
On mobility and Einstein relation for
  tracers in time-mixing random environments},  J. Stat. Phys.  {\bf 118},
   407-435   (2005).







\bibitem{liggett} T.M. Liggett, { \em Interacting Particle Systems}, Springer-Verlag,
New York (1985).


\bibitem{lu} S. Lu, { \em Hydrodynamic scaling limits with deterministic
initial configurations}, Thesis (Ph.D.)-New York University (1992).

\bibitem{n} J. Nash, {\em Continuity of solutions of parabolica and elliptic equations}, Amer. J. Math. {\bf 80}, 931-954 (1958).

\bibitem{pinsky} R. Pinsky, {\it Positive harmonic functions and diffusion}, Cambridge Univerersity Press (1995).

\bibitem{ramirez} A.F. Ram\'\i rez, {\em Relative entropy and mixing
properties of infinite dimensional diffusions}, Probab. Theory Relat. Fields
{\bf 110}, 369-395 (1998).

\bibitem{r1} A.F. Ram\'\i rez, {\em  An elementary proof of the uniqueness of invariant product
 measures for some infinite dimensional processes},
 Comptes R. Acad. Sci. {\bf 334}, S\'erie I, 1-6 (2002). 

\bibitem{r2} A.F. Ram\'\i rez, {\em Uniqueness of invariant measures for elliptic
 infinite dimensional diffusions and particle spin systems},
 ESAIM: Probabilit\'es et Statistiques {\bf 6}, 147-155 (2002).



\bibitem{varadhan} S.R. S. Varadhan. {\it Stochastic processes},
Notes based on a course given at New York University during the year 1967/68,
 Courant Institute of Mathematical Sciences, New York University, New York (1968).

\bibitem{zak} F. Zak, {\em Existence of equilibrium for Infinite System of
Interacting Diffusions}, arXiv:1502.04351 (2015).

\bibitem{z} B. Zegarlinski,
{\em Ergodicity of Markov semigroups},
Stochastic partial differential equations (Edinburgh, 1994), 312--337,
London Math. Soc. Lecture Note Sear., {\bf 216}, 
Cambridge Univ. Press, Cambridge, (1995).







\end{thebibliography}
\end{document}